\documentclass[11pt]{article}
\usepackage{amsfonts,amssymb}
\usepackage{color}
\usepackage{pstricks,pst-node}
\usepackage{pst-plot}
\catcode`\@=11
\def\marginnote#1{}

\newcount\hour
\newcount\minute
\newtoks\amorpm
\hour=\time\divide\hour by60 \minute=\time{\multiply\hour by60
\global\advance\minute by-\hour}
\edef\standardtime{{\ifnum\hour<12 \global\amorpm={am}%
        \else\global\amorpm={pm}\advance\hour by-12 \fi
        \ifnum\hour=0 \hour=12 \fi
        \number\hour:\ifnum\minute<10 0\fi\number\minute\the\amorpm}}
\edef\militarytime{\number\hour:\ifnum\minute<10 0\fi\number\minute}

\def\draftlabel#1{{\@bsphack\if@filesw {\let\thepage\relax
      \xdef\@gtempa{\write\@auxout{\string
          \newlabel{#1}{{\@currentlabel}{\thepage}}}}}\@gtempa \if@nobreak
    \ifvmode\nobreak\fi\fi\fi\@esphack} \gdef\@eqnlabel{#1}}
    \def\@eqnlabel{}
\def\@vacuum{}
\def\draftmarginnote#1{\marginpar{\raggedright\scriptsize\tt#1}}

\def\draft{
  \oddsidemargin -.5truein
  \def\@oddfoot{\footnotesize \sl preliminary draft \hfil
    \rm\thepage\hfil\sl\today\quad\militarytime}
  \let\@evenfoot\@oddfoot \overfullrule 3pt
    \let\label=\draftlabel
    \let\marginnote=\draftmarginnote
  \def\@eqnnum{(\theequation)\rlap{\kern\marginparsep\tt\@eqnlabel}%
    \global\let\@eqnlabel\@vacuum}

  }

\voffset= - 0.5in \hoffset= - 1.0in

\def\be{\begin{equation}}
\def\ee{\end{equation}}
\def\bea{\begin{eqnarray}}
\def\eea{\end{eqnarray}}
\def\<{\langle}
\def\>{\rangle}

\def\tr{{\mathrm{tr\,}}}

\def\1N{${\cal N}=1$}
\def\4N{${\cal N}=4$}

\def\CC{{\mathbb C}}

\def\ZZ{{\mathbb Z}}

\def\e{{\,\rm e}\,}

\def\bea{\begin{eqnarray}}
\def\eea{\end{eqnarray}}

\def\beq{\begin{equation}}
\def\eeq{\end{equation}}
\def\ba{\beq\begin{array}{c}}
\def\ea{\end{array}\eeq}

\@ifundefined{theorem@style}{\RequirePackage{theorem}}{}

\gdef\th@plain{\normalfont\slshape
  \def\@begintheorem##1##2{%
\item[\hskip\parindent\hskip\labelsep\theorem@headerfont ##1\ ##2\unskip.]}%
\def\@opargbegintheorem##1##2##3{%
\item[\hskip\parindent
\ifx\empty##1\else\hskip\labelsep\fi\theorem@headerfont ##1\ ##2\unskip]{\theorem@headerfont{\rm ##3}.} }}

\gdef\th@definition{\normalfont
  \def\@begintheorem##1##2{%
\item[\hskip\parindent\hskip\labelsep\theorem@headerfont ##1\ ##2\unskip.]}%
\def\@opargbegintheorem##1##2##3{%
\item[\hskip\parindent
\ifx\empty##1\else\hskip\labelsep\fi\theorem@headerfont ##1\ ##2\unskip]{\theorem@headerfont{\rm ##3}.} }}

\theoremstyle{plain}
\newtheorem{theorem}{Theorem}
\newtheorem{lemma}[theorem]{Lemma}
\newtheorem{corollary}[theorem]{Corollary}
\newtheorem{prop}[theorem]{Proposition}

\theoremstyle{definition}
\newtheorem{definition}{Definition}
\newtheorem{remark}[definition]{Remark}

\let\operatorname=\mathrm
\let\text=\mathrm

\let\wtd=\widetilde

\def\mod{\operatorname{mod}}

\def\beq{\begin{equation}}
\def\eeq{\end{equation}}
\def\bea{\begin{eqnarray}}
\def\eea{\end{eqnarray}}
\newcommand{\cpict}[3]{
\dimen1=#1\advance\dimen1 by-\hsize\divide\dimen1 by-2 \vtop to #2{
\noindent\hskip\dimen1{\special{em:graph #3.bmp}} \vfil}\hskip-2cm }

\newcommand{\tcr}{\textcolor{red}}
\newcommand{\tcb}{\textcolor{blue}}

\newcommand{\sheet}[2]{{\stackrel{{#1}}{{#2}}}}

\let\@@savethanks\thanks
\def\thanks#1{\gdef\thefootnote{\alph{footnote}}\@@savethanks{#1}}

\begin{document}

\title{Quantum ordering for quantum geodesic functions of orbifold Riemann surfaces}
\author{Leonid O. Chekhov\thanks{Steklov Mathematical Institute and  Laboratoire Poncelet,
Moscow, Russia. email: chekhov@mi.ras.ru.} \thanks{Mathematics Department, Loughborough University, UK. email: l.chekhov@lboro.ac.uk.} \,
and Marta Mazzocco\thanks{Mathematics Department, Loughborough University, Leicestershire, UK. email: m.mazzocco@lboro.ac.uk.}}

\maketitle

\begin{abstract}
We determine the explicit quantum ordering for a special class of quantum geodesic functions
corresponding to geodesics joining exactly two orbifold points or holes on a non-compact Riemann surface.
We discuss some special cases in which these quantum geodesic functions form sub--algebras of some abstract algebras defined by the
reflection equation and we extend our results to the quantisation of  matrix elements of the Fuchsian group associated to the Riemann surface in Poincar\'e uniformization. In particular we explore an interesting relation between the deformed $U_q(\mathfrak{sl}_2)$ and the Zhedanov algebra $AW(3)$.

\end{abstract}

\phantom{XXX}\hfill {\em For the 75th birthday of Sergei Petrovich Novikov}

\bigskip

\section{Introduction}
Convenient Darboux coordinates for Riemann surfaces with holes were identified in \cite{ChF1} with the shear
coordinates for an ideal triangle decomposition obtained in  \cite{Fock1} by
generalising the results obtained for punctured
Riemann surfaces in \cite{Penn1}.

These coordinates were quantized in \cite{ChF1} and, in a parametrization related to Liouville-type models,
in \cite{Kashaev}. In all cases, the quantum mapping-class group (MCG) transformations (or, the quantum flip
morphsms) that satisfy the quantum pentagon identity
involve the quantum dilogarithm function \cite{Faddeev}. Almost simultaneously, the corresponding Darboux coordinates
were identified with the $Y$-type cluster variables \cite{FZ},~\cite{FZ2}, and the cluster systems for
bordered Riemann surfaces were considered in \cite{FST}.
The above construction was generalized further to
the case of Riemann surfaces with holes and with orbifold points for the case of $\mathbb Z_2$ and $\mathbb Z_3$
orbifold points in \cite{Ch1a},~\cite{Ch2} and for the case or orbifold points of any integer order in \cite{ChSh}.

The principal object of investigation for 2D Riemann surfaces are geodesic functions: they are invariant under the
MCG action thus constituting the set of observables; in the classical case, the set of lengths of these functions
can be identified with
the spectrum of a Riemann surface, whereas a quantum Riemann surface is characterized by an algebra of quantum
geodesic functions. A simple and explicit construction of the corresponding classical geodesic functions in terms of
Darboux coordinates of decorated Teichm\"uller spaces for Riemann surfaces with holes
was proposed in \cite{ChF2}: it was shown there that all
geodesic functions are Laurent polynomials of exponentiated coordinates with positive integer coefficients; the same
remains true for Riemann surfaces with $\mathbb Z_2$ and $\mathbb Z_3$ orbifold points whereas the integrity
condition breaks in general in the case of orbifold points of arbitrary order \cite{ChSh}.

The problem of quantum ordering of a product of noncommuting operators appeared almost
simultaneously with the rise of quantum mechanics. In the context of quantum geodesic functions it
was first mentioned in \cite{ChF2} where the crucial conditions of mapping-class-group (MCG) invariance and
satisfaction of the quantum skein relations were formulated. The compatibility of these two conditions was
implicitly proved by Kashaev~\cite{Kashaev-Dehn} who constructed unitary operators of quantum
Dehn twists whose action on operators of quantum geodesic functions obviously preserves their quantum
algebra. It remained however the problem of formulating a recipe for obtaining a quantum operator in
an explicit form, likewise the Kulish, Sklyanin, and Nazarov recipe (see \cite{KulSk}, \cite{Naz}) for constructing Yangian central elements extended
to the case of reflection equation by Molev, Ragoucy, and Sorba (the quantum ordering for
twisted Yangians was constructed in \cite{MR} for the $O(n)$ case and in \cite{MRS} for the
$Sp(2n)$ case).

In the present paper, we revise the classical and quantum MCG transformations for Riemann surfaces with holes and
orbifold points and, using the results of \cite{ChP1},
construct the proper quantum ordering for a special class of geodesic functions corresponding to geodesics going
around exactly two orbifold points/holes. We use the same quantum ordering for quantizing matrix elements of the Fuchsian group
associated to the Riemann surface. In the case of {\em monodromy matrices}, Korotkin and Samtleben in~\cite{KS}
proposed an $r$-matrix structure of the Fock--Rosly type~\cite{Fock-Rosly}
which did not however satisfy Jacobi relations on monodromy matrices themselves but became
consistent on the level of  adjoint invariant elements. The quantum algebra of entries of the
matrices considered in this paper is free of this discrepansy and it is a well-defined quantum algebra. We were able to
construct these algebras for matrices in two important cases: in the case of $A_n$-algebra related to
Schlesinger systems~\cite{Dub} and in the case of the algebra related to the Painlev\'e VI equation~\cite{ChM-D4}.
In the former case, for each fixed geodesic, the matrix entries of the corresponding element in the Fuchsian group satisfy the quantum universal enveloping algebra $U_q(\mathfrak{sl}_2)$ relations, and matrix entries coming from different matrices satisfy quantum commutation relations which obey $r$-matrix type relations,  thus endowing $U_q(\mathfrak{sl}_2)\times \dots\times U_q(\mathfrak{sl}_2)$ with a well defined quantum algebra structure. We also prove that  this structure is preserved by the quantum braid group action.

In the case of the algebra related to the Painlev\'e VI equation,  for each fixed geodesic, the matrix entries of the corresponding element in the Fuchsian group satisfy  a deformed version of  the quantum universal enveloping algebra $U_q(\mathfrak{sl}_2)$ relations. This result is quite interesting as it sheds light on the relation between the quantum universal enveloping algebra $U_q(\mathfrak{sl}_2)$ and the Zhedanov algebra $AW(3)$ already explored in \cite{WZ,Ter}. Indeed, in a recent paper by M.M. \cite{MM}, it was shown that there is a natural quantisation of the monodromy group associated to the sixth Painlev\'e equation which leads to the Cherednik algebra of type $\check{C_1}C_1$. When restricted to the spherical sub-algebra the same quantization corresponds to the isomorphism between quantum geodesic functions on a Riemann sphere with four holes and Zhedanov algebra $AW(3)$  \cite{ChM-D4,IT}.

We leave the case of $D_n$-algebras related to the reflection equation with the spectral parameter to subsequent
publications.

The structure of the paper is as follows. Section~\ref{s:graph} contains a brief review of
the fat-graph description of Teichm\"uller spaces of Riemann surfaces with holes and orbifold points.
In Sec.~\ref{s:MCG}, we quantize the MCG transformations for
orbifold Riemann surfaces proposed in \cite{ChSh} thus introducing a new class of quantum MCG transformations. In
Sec.~\ref{ss:QGF}, we prove that the corresponding transformations are homogeneous for {\em matrix} products of
quantum operators. This enables us to prove in Sec.~\ref{ss:invariance} the invariance w.r.t. the quantum MCG transformations for the quantum geodesic functions corresponding to geodesics that are homeomorphic to curves
separating two orbifold points/holes from the rest of a Riemann surface. This class of geodesic functions is of
particular interest: for example, it was shown in \cite{ChM} that such geodesic functions can be identified with
elements of particular Poisson leaves of twisted Yangians (see \cite{Molev} and references therein).
Then, in Sec.~\ref{s:quantum-monodromy}, we construct quantum  matrices on rooted fat graphs
that are MCG-invariant w.r.t. all transformations that leave the root edge intact. The quantum algebra for
the $A_n$ case is constructed in Sec.~\ref{ss:An-monodromy} and we represent it in
the $R$-matrix form in Sec.~\ref{ss:R-matrix} thus constructing the consistent quantum algebra for
monodromy matrices. We present the action of the quantum braid group
transformations on the level of monodromy matrices in Sec.~\ref{ss:braid}.
In Sec.~\ref{ss:PVI-monodromy}, we find the quantum
algebra for monodromy matrices of the Painlev\'e VI equation.

\section{Combinatorial description of ${\mathfrak T}^{H}_{g,s,r}$}\label{s:graph}

\subsection{Fat graph description for Riemann surfaces with holes and ${\mathbb Z}_p$ orbifold points}

\begin{definition}\label{def-pend}
We call a fat graph (a graph with the prescribed cyclic ordering of edges
entering each vertex) $\Gamma_{g,s,r}$ a {\em spine of the Riemann surface} $\Sigma_{g,s,r}$
with $g$ handles, $s>0$ holes, and $r$ orbifold points of the corresponding orders $p_i$, $i=1,\dots,r$, if
\begin{itemize}
\item[(a)] this graph can be embedded  without self-intersections in $\Sigma_{g,s,r}$;
\item[(b)] all vertices of $\Gamma_{g,s,r}$ are three-valent except exactly $r$
one-valent vertices (endpoints of ``pending'' edges), which are placed at the corresponding
orbifold points;
\item[(c)] upon cutting along all edges of $\Gamma_{g,s,r}$ the Riemann surface
$\Sigma_{g,s,r}$ splits into $s$ polygons each containing exactly one hole and being
simply connected upon contracting  this hole.
\end{itemize}
Because every pending edge ``protrudes'' towards the interior of some face of the graph and every face contains exactly one hole, the
above fat graph determines a natural partition of the set of orbifold points into nonintersecting (maybe empty)
subsets $\delta_k$, $k=1,\dots,s$ of orbifold points incident to the corresponding face (boundary component, or hole). 
Edges of the above graph are labeled by distinct integers $\alpha=1,2,\dots,6g-6+3s+2r$, and we set
a real number $Z_\alpha$ into correspondence to the $\alpha$th edge.
\end{definition}

The first homotopy groups $\pi_1(\Sigma_{g,s,r})$ and $\pi_1(\Gamma_{g,s,r})$ coincide because
each closed path in $\Sigma_{g,s,r}$ can be homotopically transformed to a closed path in $\Gamma_{g,s,r}$
(taking into account paths that go around orbifold points)
in a unique way. The standard statement in hyperbolic geometry is that conjugacy classes of elements of
a Fuchsian group $\Delta_{g,s,r}$ are in the 1-1 correspondence with homotopy
classes of closed paths in the Riemann surface $\Sigma_{g,s,r}={\mathbb H}^2_+/\Delta_{g,s,r}$ and that the
``actual'' length $\ell_\gamma$
of a hyperbolic element $\gamma\in\Delta_{g,s,r}$ coincides with the minimum length of
curves from the corresponding homotopy class; it is then the length of a unique closed
geodesic line belonging to this class.

The real numbers $Z_\alpha$ in Definition~\ref{def-pend} are the $h$-lengths (logarithms of cross-ratios)
\cite{Penn1}: they are called the {\em (Thurston) shear
coordinates} \cite{ThSh},\cite{Bon2} in the case of punctured Riemann surface. Below we identify these
shear coordinates with coordinates of the decorated Teichm\"uller space ${\mathfrak T}^{H}_{g,s,r}$.

\subsection{The Fuchsian group $\Delta_{g,s,r}$ and geodesic functions}\label{ss:geodesic}

We now describe combinatorially the conjugacy classes of the Fuchsian group $\Delta_{g,s,r}$.
Every time the path homeomorphic to a (closed) geodesic $\gamma$ passes along the edge with the label $\alpha$ we
insert~\cite{Fock1} the so-called {\it edge matrix}, i.e. the matrix of M\"obius transformation:
\be
\label{XZ} X_{Z_\alpha}=\left(
\begin{array}{cc} 0 & -\e^{Z_\alpha/2}\\
                \e^{-Z_\alpha/2} & 0\end{array}\right)
\ee
into the corresponding string of matrices. We also have the ``right'' and ``left'' turn matrices
to be set in proper places when a path makes corresponding turns at three-valent vertices,
\be
\label{R}
R=\left(\begin{array}{cc} 1 & 1\\ -1 & 0\end{array}\right), \qquad
L= R^2=\left(\begin{array}{cc} 0 & 1\\ -1 &
-1\end{array}\right).
\ee

\begin{figure}[tb]
{\psset{unit=0.7}
\begin{pspicture}(-4,-3)(2,2)
\pcline[linewidth=1pt](-1.5,1)(1.5,1)
\pcline[linewidth=1pt](-1.5,0)(-0.5,0)
\pcline[linewidth=1pt](1.5,0)(0.5,0)
\pcline[linewidth=1pt](-.5,0)(-.5,-1)
\pcline[linewidth=1pt](.5,0)(.5,-1)
\rput(0,-1){\makebox(0,0){$\bullet$}}
\pcline[linecolor=blue, linestyle=dashed, linewidth=1pt]{<-}(-1.5,0.2)(-.5,0.2)
\psarc[linecolor=blue, linestyle=dashed, linewidth=1pt](-.5,0){.2}{0}{90}
\pcline[linecolor=blue, linestyle=dashed, linewidth=1pt](-.3,0)(-.3,-1)
\psarc[linecolor=blue, linestyle=dashed, linewidth=1pt](0,-1){.3}{-180}{0}
\pcline[linecolor=blue, linestyle=dashed, linewidth=1pt](.3,0)(.3,-1)
\psarc[linecolor=blue, linestyle=dashed, linewidth=1pt](.5,0){.2}{90}{180}
\pcline[linecolor=blue, linestyle=dashed, linewidth=1pt]{<-}(.5,0.2)(1.5,0.2)
\rput(-1.3,1.2){\makebox(0,0)[cb]{$X$}}
\rput(1.3,1.2){\makebox(0,0)[cb]{$Y$}}
\rput(0.8,-.7){\makebox(0,0){$Z$}}
\rput(0,-1.8){\makebox(0,0){$F_\omega$}}
\rput(0,-2.4){\makebox(0,0){(a)}}
\end{pspicture}
\begin{pspicture}(-2,-3)(2,2)
\pcline[linewidth=1pt](-1.5,1)(1.5,1)
\pcline[linewidth=1pt](-1.5,0)(-0.5,0)
\pcline[linewidth=1pt](1.5,0)(0.5,0)
\pcline[linewidth=1pt](-.5,0)(-.5,-1)
\pcline[linewidth=1pt](.5,0)(.5,-1)
\rput(0,-1){\makebox(0,0){$\bullet$}}
\pcline[linecolor=blue, linestyle=dashed, linewidth=1pt]{<-}(-1.5,0.2)(-.5,0.2)
\psarc[linecolor=blue, linestyle=dashed, linewidth=1pt](-.5,0.05){.15}{0}{90}
\pcline[linecolor=blue, linestyle=dashed, linewidth=1pt](-.35,0.05)(-.35,-1)
\psarc[linecolor=blue, linewidth=0.7pt](0.075,-1){0.275}{-180}{0}
\psarc[linecolor=blue, linewidth=0.7pt](-0.075,-1){0.275}{-180}{0}
\psarc[linecolor=blue, linewidth=0.7pt](0,-1){0.2}{0}{180}
\pcline[linecolor=blue, linestyle=dashed, linewidth=1pt](.35,0.05)(.35,-1)
\psarc[linecolor=blue, linestyle=dashed, linewidth=1pt](-.4,0.05){.75}{0}{90}
\pcline[linecolor=blue, linestyle=dashed, linewidth=1pt]{->}(-1.5,0.8)(-.4,0.8)
\rput(-1.3,1.2){\makebox(0,0)[cb]{$X$}}
\rput(1.3,1.2){\makebox(0,0)[cb]{$Y$}}
\rput(0.8,-.7){\makebox(0,0){$Z$}}
\rput(0,-1.8){\makebox(0,0){$-F^2_\omega$}}
\rput(0,-2.4){\makebox(0,0){(b)}}
\end{pspicture}
\begin{pspicture}(-2,-3)(2,2)
\pcline[linewidth=1pt](-1.5,1)(1.5,1)
\pcline[linewidth=1pt](-1.5,0)(-0.5,0)
\pcline[linewidth=1pt](1.5,0)(0.5,0)
\pcline[linewidth=1pt](-.5,0)(-.5,-1)
\pcline[linewidth=1pt](.5,0)(.5,-1)
\rput(0,-1){\makebox(0,0){$\bullet$}}
\pcline[linecolor=blue, linestyle=dashed, linewidth=1pt]{->}(1.5,0.2)(.5,0.2)
\psarc[linecolor=blue, linestyle=dashed, linewidth=1pt](.5,0.1){.1}{90}{180}
\pcline[linecolor=blue, linestyle=dashed, linewidth=1pt](-.4,0.1)(-.4,-1)
\psarc[linecolor=blue, linewidth=0.7pt](0.05,-1){.35}{-180}{0}
\psarc[linecolor=blue, linewidth=0.7pt](-0.05,-1){.35}{-180}{0}
\psarc[linecolor=blue, linewidth=0.7pt](0,-1){.2}{-180}{0}
\psarc[linecolor=blue, linewidth=0.7pt](0.05,-1){.25}{0}{180}
\psarc[linecolor=blue, linewidth=0.7pt](-0.05,-1){.25}{0}{180}
\pcline[linecolor=blue, linestyle=dashed, linewidth=1pt](.4,0.1)(.4,-1)
\psarc[linecolor=blue, linestyle=dashed, linewidth=1pt](.3,0.1){.7}{90}{180}
\pcline[linecolor=blue, linestyle=dashed, linewidth=1pt]{<-}(1.5,0.8)(.3,0.8)
\rput(-1.3,1.2){\makebox(0,0)[cb]{$X$}}
\rput(1.3,1.2){\makebox(0,0)[cb]{$Y$}}
\rput(0.8,-.7){\makebox(0,0){$Z$}}
\rput(0,-1.8){\makebox(0,0){$F^3_\omega$}}
\rput(0,-2.4){\makebox(0,0){(c)}}
\end{pspicture}
}
\caption{\small Part of a graph with a pending edge.
Its endpoint with the orbifold point is directed toward the interior
of the boundary component this point is associated with.
The variable $Z$ corresponds to the respective pending edge. We present
three typical examples of geodesics undergoing single (a), double (b),
and triple (c) rotations at the ${\mathbb Z}_p$ orbifold point.}
\label{fi:corner}
\end{figure}
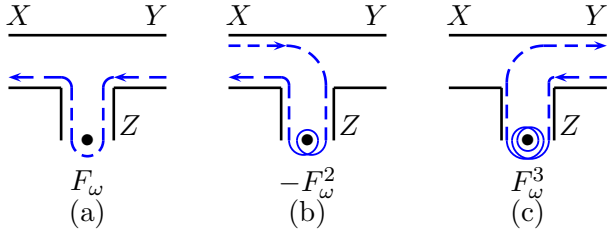

When orbifold points are present, the Fuchsian group contains besides hyperbolic elements also elliptic
elements corresponding to rotations about these orbifold points.
The corresponding generators ${\wtd F}_i$, $i=1,\dots,r$, of the rotations through $2\pi/p_i$
are conjugates of the matrices
\be
\label{F-p}
{\wtd F}_i=U_iF_{\omega_i}U_i^{-1},\qquad F_\omega:=
\left(\begin{array}{cc} 0 & 1\\ -1 & -w\end{array}\right),\quad w=2\cos{\pi/p}\ \hbox{for some}\ p\ge 2.
\ee
New elements of the
Fuchsian group correspond to rotations of geodesics when going around orbifold points
indicated by dot-vertices;
for a ${\mathbb Z}_p$ orbifold point we then insert the above matrix $F_\omega$
into the corresponding string of $2\times2$-matrices (when we go around the orbifold point counterclockwise
as in Fig.~\ref{fi:corner}(a)). When going around an orbifold point
$k$ times we insert the matrix $(-1)^{k+1}F_\omega^k$ into the product of $2\times2$-matrices.
For example, parts of geodesic functions
in the three cases in Fig.~\ref{fi:corner} read
\be
\label{XZFXZ}
\begin{array}{ll}
\hbox{(a)}\quad & \dots  X_XLX_ZF_\omega X_ZLX_Y\dots, \\
\hbox{(b)}\quad & \dots  X_XLX_Z(-F^2_\omega)X_ZRX_X\dots, \\
\hbox{(c)}\quad & \dots  X_YRX_Z(F^3_\omega)X_ZLX_Y\dots. \\
\end{array}
\ee
Note that for $\omega_p=2\cos{\pi/p}$, $F_{\omega_p}^p=(-1)^{p-1}{\mathbb E}$,
so going around the ${\mathbb Z}_p$ orbifold point $p$ times merely corresponds to avoiding this orbifold point
due to the simple equality (note that $X_S^2=-{\mathbb E}$ and $L^2=-R$)
$$
X_XLX_Z(-1)^{p-1}F_{\omega_p}^pX_ZLX_Y=X_XLX_Z^2LX_Y=-X_XL^2X_Y=X_XRX_Y.
$$
(For the ${\mathbb Z}_2$ orbifold points
this pattern was first proposed by Fock and Goncharov \cite{FG}; the graph morphisms were described in \cite{Ch2}.)

An element of a Fuchsian group has then the typical structure
\be
\label{Pgamma}
P_{\gamma}=LX_{Y_n}RX_{Y_{n-1}}\cdots RX_{Y_2}LX_{Z_1}(-1)^{k+1}F^k_{\omega_i} X_{Z_1}RX_{Y_1},
\ee
where $Y_i$ are variables of ``internal'' edges and $Z_j$ are those of pending edges.
In the corresponding {\em geodesic function}
\be
\label{G}
G_{\gamma}\equiv \tr P_\gamma=2\cosh(\ell_\gamma/2),
\ee
$\ell_\gamma$ is the actual length of the closed
geodesic on the Riemann surface.

\begin{remark}\label{rm-positivity}
Note that the combinations
\bea
RX_Y&=&\left(\begin{array}{cc} e^{-Y/2} & -e^{Y/2}\\ 0 & e^{Y/2}\end{array}\right),
\quad \hbox{and}\quad
LX_Y=\left(\begin{array}{cc} e^{-Y/2} & 0 \\ -e^{-Y/2} & e^{Y/2}\end{array}\right),\nonumber\\
RX_ZF_\omega X_Z&=&\left(\begin{array}{cc} e^{-Z}+\omega & -e^{Z}\\ -\omega & e^{Z}\end{array}\right)
\quad\hbox{and}\quad
LX_ZF_\omega X_Z=\left(\begin{array}{cc} e^{-Z} & 0\\ -e^{-Z}-\omega & e^{Z}\end{array}\right),
\nonumber
\eea
as well as products of any number of these matrices have the sign structure
$\left(\begin{array}{cc} + & -\\ - & +\end{array}\right)$, so the trace of any of $P_\gamma$ with first powers of
$F_\omega$ is a
sum of exponentials with positive integer coefficients; this sum always include the terms
$$
e^{Y_1/2+\cdots+Y_n/2+Z_1+\cdots +Z_s}\quad \hbox{and}\quad e^{-Y_1/2-\cdots-Y_n/2-Z_1-\cdots -Z_s}
$$
being therefore always greater or equal
two thus describing a hyperbolic or parabolic element;
the latter is possible only if
$Y_1+\cdots+Y_n+Z_1+\cdots +Z_s=0$ and only if all the turn matrices in (\ref{Pgamma}) are left-turn ones
(for our choice of $F_\omega$; when going in the opposite direction we must insert $F_\omega^{-1}$, and
all the turn matrices must be right-turn ones),
which corresponds to a path going along the boundary of a face; all such paths are homeomorphic
to the hole boundaries, and the condition that the sum of $Y_i$ and $Z_j$ equals zero along this path
indicates the degeneration of a hole into a puncture. Fot the complete proof of hyperbolicity of
all elements not homeomorphic to going several times around a single orbifold point, see \cite{ChSh}.
\end{remark}

The group generated by elements (\ref{F-p}) together with translations along $A$- and $B$-cycles
and around holes not necessarily produces a regular (metrizable) surface because its action
is not necessarily discrete. The necessary and
sufficient conditions for producing a {\em regular} surface
in terms of graphs were formulated as a theorem in \cite{ChSh}. (We call a Riemann surface
regular if it is locally a smooth constant-curvature surface everywhere except exactly $r$ orbifold points.)
Essentially, it was proved that we obtain a regular Riemann surface for any choice of real
numbers $Y_i$ and $Z_\alpha$
from Definition~\ref{def-pend} and vice versa, for any regular Riemann surface $\Sigma_{g,s,r}={\mathbb H}^2_+/\Delta_{g,s,r}$ and for any fat graph that is a spine of this surface we find a set of real numbers
$Y_i$ and $Z_\alpha$ that are the shear coordinates for an ideal triangle decomposition of the above surface. This set is not unique; however the equivalent sets are related by discrete modular group action, and we can therefore
identify the $(6g-6+3s+2r)$-tuple of real coordinates $\{Y_i,Z_\alpha\}$
with the coordinates of the decorated
Teichm\"uller space ${\mathfrak T}^{H}_{g,s,r}$ (the decoration assigns positive or negative signs to
every hole with nonzero perimeter). The
lengths of geodesics on $\Sigma_{g,s,r}$ are given by traces of products (\ref{Pgamma})
corresponding to paths in the corresponding spine.

We see that every spine $\Gamma_{g,s,r}$ provides a parameterization of the (decorated) Teichm\"uller space ${\mathfrak T}^H_{g,s,r}$. Transitions between different parameterizations are formulated in terms of flips
(mutations) of edges: any two spines from the given topological class are related by a finite sequence of flips.
We therefore identify flips of edges with the action of the MCG to be constructed in the next
section.

\section{Quantum mapping class group transformations for surfaces with orbifold points}\label{s:MCG}

\subsection{Poisson structure}\label{ss:Poisson}

One of the most attractive properties of the graph description is a very simple Poisson algebra on the set
of parameters $Z_\alpha$, $\alpha=1,\dots, 6g-6+3s+2r$
(we do not distinguish here the parameters of internal and pending edges).

\begin{theorem}\label{th-WP} In the coordinates $Z_\alpha $ on any fixed spine
corresponding to a surface with or without orbifold points,
the Weil--Petersson bracket $B_{{\mbox{\tiny WP}}}$ reads
\be
\label{WP-PB}
\bigl\{f({\mathbf Z}),g({\mathbf Z})\bigr\}=\sum_{{\hbox{\small 3-valent} \atop \hbox{\small vertices $\alpha=1$} }}^{4g+2s+|\delta|-4}
\,\sum_{i=1}^{3 \mod 3}
\left(\frac{\partial f}{\partial Z_{\alpha_i}} \frac{\partial g}{\partial Z_{\alpha_{i+1}}}
- \frac{\partial g}{\partial Z_{\alpha_i}} \frac{\partial f}{\partial Z_{\alpha_{i+1}}}\right),
\ee
where the sum ranges all the three-valent vertices of a graph and
$\alpha_i$ are the labels of the cyclically (clockwise)
ordered ($\alpha_4\equiv \alpha_1 $) edges incident to the vertex
with the label $\alpha$. The formula (\ref{WP-PB}) is insensitive to whether these edges are internal or
pending edges. This bracket gives rise to the {\em Goldman
bracket} on the space of geodesic length functions \cite{Gold}.
\end{theorem}

The center of the Poisson algebra {\rm(\ref{WP-PB})} is generated by
elements of the form $\sum Z_\alpha$, where the sum ranges all edges
of $\Gamma_{g,s,r} $ belonging to the same boundary component
taken with multiplicities.
This means, in particular, that each pending edge, irrespectively on the type of orbifold point it corresponds to,
contributes twice to such sums. The dimension of this center is obviously $s$.

Note that for the path homeomorphic to
the hole boundary, for any number of insertions of matrices $F_{p_i}$ with any $p_i$, we have
\bea
&{}&\hbox{tr\,}\bigl[LX_{Z_1}LX_{Z_2}\cdots LX_{Z_k}F_{\omega}X_{Z_k}L\cdots LX_{Z_{n-1}}LX_{Z_n}\bigr]
\nonumber\\
&{}&\quad = 2\cosh \bigl[Z_1/{2}+{Z_2}/{2}+\dots +Z_k+\dots +{Z_{n-1}}/{2}+{Z_n}/{2}\bigr].
\nonumber
\eea

\subsection{Flip morphisms of fat graphs}\label{ss:flip}

We now describe all mapping class group transformations
that enable us to change numbers $|\delta_k|$ of orbifold points associated with the $k$th
hole, change the cyclic ordering inside any of the sets $\delta_k$, flip any inner edge
of the graph and, eventually, change the decoration, i.e., orientation of the geodesic
spiraling to the hole perimeter (in the case where we have more than one hole).
We can therefore establish a morphism between any two of the graphs belonging to the same class $\Gamma_{g,s,r}$
with the same (unordered) sets of orbifold point orders $\{p_i\}_{i=1}^r$.

\subsubsection{Flipping inner edges}\label{sss:mcg}

\begin{figure}[tb]
\begin{pspicture}(-6,-3)(4,3){
\newcommand{\FLIP}{%
{\psset{unit=1}
\psline[linewidth=18pt,linecolor=blue](0,-1)(0,1)
\psline[linewidth=18pt,linecolor=blue](0,1)(1.5,2)
\psline[linewidth=18pt,linecolor=blue](0,1)(-1.5,2)
\psline[linewidth=18pt,linecolor=blue](0,-1)(1.5,-2)
\psline[linewidth=18pt,linecolor=blue](0,-1)(-1.5,-2)
\psline[linewidth=14pt,linecolor=white](0,-1)(0,1)
\psline[linewidth=14pt,linecolor=white](0,1)(1.5,2)
\psline[linewidth=14pt,linecolor=white](0,1)(-1.5,2)
\psline[linewidth=14pt,linecolor=white](0,-1)(1.5,-2)
\psline[linewidth=14pt,linecolor=white](0,-1)(-1.5,-2)
}
}
\rput(-2.5,0){\FLIP}
\psline[linewidth=2pt]{<->}(-0.5,0)(0.5,0)
\rput{90}(3.5,0){\FLIP}
\rput(-2.5,0){
\rput(-1.1,2.2){\makebox(0,0)[lb]{$A$}}
\rput(1.1,2.2){\makebox(0,0)[rb]{$B$}}
\rput(-0.8,0){\makebox(0,0)[lc]{$Z$}}
\rput(1.1,-2.2){\makebox(0,0)[rt]{$C$}}
\rput(-1.1,-2.2){\makebox(0,0)[lt]{$D$}}
\psline[linecolor=red](-1.5,-2)(-.15,-1.1)
\psline[linecolor=red](1.5,2)(.15,1.1)
\psbezier[linecolor=red](-.15,-1.1)(0.15,-0.9)(-0.15,0.9)(.15,1.1)
\rput(-.1,.15){\psline[linecolor=red](-1.5,-2)(-.15,-1.1)}
\rput(-.1,-.15){\psline[linecolor=red](-1.5,2)(-.15,1.1)}
\psbezier[linecolor=red](-.25,-.95)(-0.1,-0.85)(-0.1,0.85)(-.25,.95)
\rput(.08,-.12){\psline[linecolor=red](-1.5,-2)(-.15,-1.1)}
\rput(-.08,-.12){\psline[linecolor=red](1.5,-2)(.15,-1.1)}
\psbezier[linecolor=red](-.07,-1.22)(0,-1.176)(0,-1.176)(.07,-1.22)
}
\rput(3.5,0){
\rput(-2.2,-2){\makebox(0,0)[lt]{$D - \phi(-Z)$}}
\rput(2.2,-2){\makebox(0,0)[rt]{$C+\phi(Z)$}}
\rput(2.2,2){\makebox(0,0)[rb]{$B-\phi(-Z)$}}
\rput(-2.2,2){\makebox(0,0)[lb]{$A+\phi(Z)$}}
\rput(0,0.5){\makebox(0,0)[cb]{$-Z$}}
}
\rput{90}(3.5,0){
\psline[linecolor=red](1.5,-2)(.15,-1.1)
\psline[linecolor=red](-1.5,2)(-.15,1.1)
\psbezier[linecolor=red](.15,-1.1)(-0.15,-0.9)(0.15,0.9)(-.15,1.1)
\rput(-.1,.15){\psline[linecolor=red](-1.5,-2)(-.15,-1.1)}
\rput(-.1,-.15){\psline[linecolor=red](-1.5,2)(-.15,1.1)}
\psbezier[linecolor=red](-.25,-.95)(-0.1,-0.85)(-0.1,0.85)(-.25,.95)
}
\rput{270}(3.5,0){
\rput(.08,-.12){\psline[linecolor=red](-1.5,-2)(-.15,-1.1)}
\rput(-.08,-.12){\psline[linecolor=red](1.5,-2)(.15,-1.1)}
\psbezier[linecolor=red](-.07,-1.22)(0,-1.176)(0,-1.176)(.07,-1.22)
}
\rput(-2.5,0){
\put(-1.8,1.9){\makebox(0,0)[cc]{\hbox{\tcr{\small$1$}}}}
\put(-1.8,1.9){\pscircle[linecolor=red]{.2}}
\put(1.7,2.1){\makebox(0,0)[cc]{\hbox{\tcr{\small$2$}}}}
\put(1.7,2.1){\pscircle[linecolor=red]{.2}}
\put(1.6,-2.3){\makebox(0,0)[cc]{\hbox{\tcr{\small$3$}}}}
\put(1.6,-2.3){\pscircle[linecolor=red]{.2}}
}
\rput(3.5,0){
\put(-2.2,1.7){\makebox(0,0)[cc]{\hbox{\tcr{\small$1$}}}}
\put(-2.2,1.7){\pscircle[linecolor=red]{.2}}
\put(2.2,1.7){\makebox(0,0)[cc]{\hbox{\tcr{\small$2$}}}}
\put(2.2,1.7){\pscircle[linecolor=red]{.2}}
\put(2.1,-1.6){\makebox(0,0)[cc]{\hbox{\tcr{\small$3$}}}}
\put(2.1,-1.6){\pscircle[linecolor=red]{.2}}
}
}
\end{pspicture}
\caption{\small Flip on the shear coordinates $Z_\alpha$. The outer edges can be
pending, but the edge undergoing the flip must be an internal
edge with distinct endpoints. We also indicate the correspondences between geodesic paths under the flip.}
\label{fi:flip}
\end{figure}
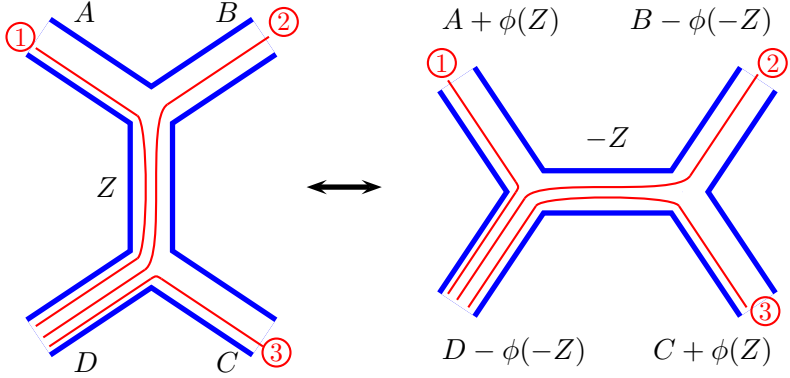

Given a spine $\Gamma$ of $\Sigma$ and
assuming that the internal
edge $\alpha$ has distinct endpoints, we may produce
another spine $\Gamma _\alpha$ of $\Sigma$ by contracting and expanding edge $\alpha$ of
$\Gamma $, the edge labeled $Z$ in Figure~\ref{fi:flip}.
We say that $\Gamma _\alpha$ arises from $\Gamma$ by a
{\it Whitehead move} (or flip) along the edge $\alpha$.
A labeling of edges of the spine $\Gamma$ implies a natural labeling of edges of the
spine $\Gamma_\alpha$; we then obtain a morphism between the spines $\Gamma$ and $\Gamma_\alpha$.

It was shown in \cite{ChF1} that setting $\phi (Z)={\rm log}(1+e^Z)$ and adopting the notation of Fig.~\ref{fi:flip}
for shear coordinates of nearby edges, the effect of a flip is
\bea
W_Z\,:\ (A,B,C,D,Z)&\to& (A+\phi(Z), B-\phi(-Z), C+\phi(Z), D-\phi(-Z), -Z)\nonumber\\
&:=&({\tilde A},{\tilde B},{\tilde C},{\tilde D},{\tilde Z}).
\label{abc}
\eea
In various cases where the edges are not distinct, we have:
if $A=C$, then $\tilde A=A+2\phi(Z)$;
if $B=D$, then $\tilde B=B-2\phi(-Z)$;
if $A=B$ (or $C=D$), then $\tilde A=A+Z$ (or $\tilde C=C+Z$);
if $A=D$ (or $B=C$), then $\tilde A=A+Z$ (or $\tilde B=B+Z$).
Any subset of edges $A$, $B$, $C$, and $D$ can be pending edges of the graph.

We have the lemma establishing the properties of
invariance w.r.t. the flip morphisms~\cite{ChF2}.

\begin{lemma} \label{lem-abc}
Transformation~{\rm(\ref{abc})} preserves
the traces of products over paths {\rm(\ref{G})} (the geodesic functions) and
transformation~{\rm(\ref{abc})} simultaneously preserves
Poisson structure {\rm(\ref{WP-PB})} on the shear coordinates.
\end{lemma}

The proof of this lemma is based on the following useful {\em matrix} equalities, which are valid
without taking the trace (they correspond to three geodesic cases in Fig.~\ref{fi:flip}):
\bea
X_DRX_ZRX_{A}&=&X_{\tilde A}RX_{\tilde D},\label{mutation1}\\
X_DRX_ZLX_B&=&X_{\tilde D}LX_{\tilde Z}RX_{\tilde B},\label{mutation2}\\
X_CLX_D&=&X_{\tilde C}LX_{\tilde Z}LX_{\tilde D}.\label{mutation3}
\eea

\subsubsection{Flipping pending edges}\label{sss:pending}

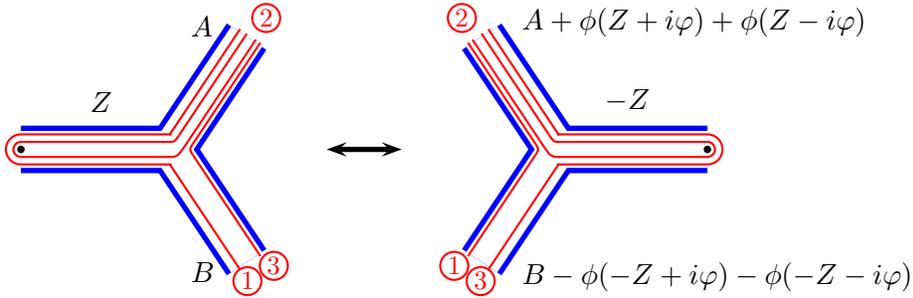
\begin{figure}[tb]
\begin{pspicture}(-6,-3)(4,3){
\newcommand{\FLIP}{%
{\psset{unit=1}
\psline[linewidth=18pt,linecolor=blue](-2,0)(0,0)
\psline[linewidth=18pt,linecolor=blue](0,0)(1,1.5)
\psline[linewidth=18pt,linecolor=blue](0,0)(1,-1.5)
\psline[linewidth=14pt,linecolor=white](-2,0)(0,0)
\psline[linewidth=14pt,linecolor=white](0,0)(1,1.5)
\psline[linewidth=14pt,linecolor=white](0,0)(1,-1.5)
\rput(-2,0){\pscircle*{0.05}}
\psarc[linecolor=red](-2,0){0.1}{90}{270}
\psarc[linecolor=red](-2,0){0.2}{90}{270}
\psline[linecolor=red](-2,0.1)(-0.01,0.1)
\psline[linecolor=red](-2,0.2)(-0,0.2)
\psline[linecolor=red](-2,-0.1)(-0.01,-0.1)
\psline[linecolor=red](-2,-0.2)(-0,-0.2)
}
}
\rput(-2.5,0){\FLIP}
\psline[linewidth=2pt]{<->}(-0.5,0)(0.5,0)
\rput{180}(2.5,0){\FLIP}
\rput(-2.5,0){
\rput(0.2,1.5){\makebox(0,0)[lb]{$A$}}
\rput(-1,0.5){\makebox(0,0)[cb]{$Z$}}
\rput(0.2,-1.5){\makebox(0,0)[lt]{$B$}}
}
\rput(3.5,0){
\rput(-1.4,-1.5){\makebox(0,0)[lt]{$B - \phi(-Z+i\varphi)-\phi(-Z-i\varphi)$}}
\rput(-1.4,1.5){\makebox(0,0)[lb]{$A+\phi(Z+i\varphi)+\phi(Z-i\varphi)$}}
\rput(0,0.5){\makebox(0,0)[cb]{$-Z$}}
}
\rput{90}(3.5,0){
\rput(-.1,-.15){\psline[linecolor=red](-1.5,2)(-.09,1.06)}
\rput(.1,-.15){\psline[linecolor=red](1.5,2)(.09,1.06)}
}
\rput{270}(-3.5,0){
\rput(-.1,-.15){\psline[linecolor=red](-1.5,2)(-.09,1.06)}
\rput(.1,-.15){\psline[linecolor=red](1.5,2)(.09,1.06)}
}
\rput{270}(3.5,0){
\rput(.08,-.12){\psline[linecolor=red](-1.5,-2)(-.15,-1.1)}
\rput(-.08,-.12){\psline[linecolor=red](1.5,-2)(.15,-1.1)}
\psbezier[linecolor=red](-.07,-1.22)(0,-1.176)(0,-1.176)(.07,-1.22)
}
\rput{90}(-3.5,0){
\rput(.08,-.12){\psline[linecolor=red](-1.5,-2)(-.15,-1.1)}
\rput(-.08,-.12){\psline[linecolor=red](1.5,-2)(.15,-1.1)}
\psbezier[linecolor=red](-.07,-1.22)(0,-1.176)(0,-1.176)(.07,-1.22)
}
\rput(-2.5,0){
\rput(-.06,.04){\psline[linecolor=red](1,1.5)(0.1,.15)}
\rput(.06,-.04){\psline[linecolor=red](1,1.5)(0.1,.15)}
\psbezier[linecolor=red](.04,.19)(-0.02,0.1)(-0.02,0.1)(-0.1,0.1)
\psbezier[linecolor=red](.16,.11)(0.02,-0.1)(0.02,-.1)(-0.1,-.1)
}
\rput(2.5,0){
\rput(.06,.04){\psline[linecolor=red](-1,1.5)(-0.1,.15)}
\rput(-.06,-.04){\psline[linecolor=red](-1,1.5)(-0.1,.15)}
\psbezier[linecolor=red](-.04,.19)(0.02,0.1)(0.02,0.1)(0.1,0.1)
\psbezier[linecolor=red](-.16,.11)(-0.02,-0.1)(-0.02,-.1)(0.1,-.1)
}
\rput(-3.5,0){
\put(2.3,-1.6){\makebox(0,0)[cc]{\hbox{\tcr{\small$3$}}}}
\put(2.3,-1.6){\pscircle[linecolor=red]{.2}}
\put(2.2,1.7){\makebox(0,0)[cc]{\hbox{\tcr{\small$2$}}}}
\put(2.2,1.7){\pscircle[linecolor=red]{.2}}
\put(1.95,-1.8){\makebox(0,0)[cc]{\hbox{\tcr{\small$1$}}}}
\put(1.95,-1.8){\pscircle[linecolor=red]{.2}}
}
\rput(3.5,0){
\put(-2.3,-1.6){\makebox(0,0)[cc]{\hbox{\tcr{\small$1$}}}}
\put(-2.3,-1.6){\pscircle[linecolor=red]{.2}}
\put(-2.2,1.7){\makebox(0,0)[cc]{\hbox{\tcr{\small$2$}}}}
\put(-2.2,1.7){\pscircle[linecolor=red]{.2}}
\put(-1.95,-1.8){\makebox(0,0)[cc]{\hbox{\tcr{\small$3$}}}}
\put(-1.95,-1.8){\pscircle[linecolor=red]{.2}}
}
}
\end{pspicture}
\caption{\small
The transformation of dual variables ($h$-lengths) when flipping a pending edge;
$w=2\cos(\pi/p)$ and $\varphi=\pi/p$. We indicate the orbifold point by the bullet and demonstrate the
changing of geodesic functions upon flipping the edge.
}
\label{fi:interchange-p-dual}
\end{figure}

\begin{lemma} \label{lem-pending1}
The transformation~in Fig.~\ref{fi:interchange-p-dual}
\be
\label{morphism-pending}
\{\tilde A,\tilde B,\tilde Z\}:= \{A+\phi(Z+i\varphi)+\phi(Z-i\varphi), B-\phi(-Z+i\varphi)-\phi(-Z-i\varphi),-Z\},
\ee
where $\phi(x)=\log(1+e^x)$ and $w=e^{i\varphi}+e^{-i\varphi}$,
is the morphism of the space
${\mathfrak T}_{g,s,r}^H$. These
morphisms preserve both Poisson structures {\rm(\ref{WP-PB})} and the geodesic
functions. In Fig.~\ref{fi:interchange-p-dual} any (or both) of variables $A$ and $B$ can be
variables of pending edges (the transformation formula is insensitive to it).
\end{lemma}

In this case, again, verifying the preservation of Poisson relations (\ref{WP-PB}) is simple, whereas
for traces over paths we have three different cases, and in each of these cases we again obtain
$2\times2$-{\em matrix} equalities to be verified directly.\footnote{In \cite{ChSh}, the powers of the matrix $F_p$ were incorrect in the corresponding formula.} We let
$$
\Omega=\left(
         \begin{array}{cc}
           a & c \\
           -c & a-wc \\
         \end{array}
       \right),\quad \forall a,c\in \CC,
$$
denote any matrix commuting with $F_p$ (in particular, $F_p^k$, $k\in{\ZZ}$). We then have the following
matrix equalities:
\bea
X_{A}LX_Z(F_p\Omega) X_ZLX_{B}&=&
X_{{\tilde A}}RX_{\tilde Z}(-\Omega) X_{\tilde Z}RX_{{\tilde B}},
\label{pend-mutation1}\\
X_{A}LX_Z\Omega X_ZRX_{A}&=&
X_{{\tilde A}}RX_{\tilde Z}(-\Omega)X_{\tilde Z}LX_{{\tilde A}},
\label{pend-mutation2}\\
X_{B}RX_Z\Omega X_ZLX_{B}&=&
X_{{\tilde B}}LX_{\tilde Z}(-\Omega)X_{\tilde Z}RX_{{\tilde B}}.
\label{pend-mutation3}
\eea

Using flip morphisms in Fig.~\ref{fi:interchange-p-dual} and in formula
(\ref{abc}), we establish a morphism between any two algebras
corresponding to surfaces of the same genus, same number of boundary
components, and same numbers of ${\mathbb Z}_p$-orbifold
points of each sort $p$; the distribution of latter into the boundary components as well as
the cyclic ordering inside each of the boundary component can be arbitrary.

It is a standard tool that if, after a
series of morphisms, we come to a graph of the same combinatorial
type as the initial one (disregarding labeling of edges but distinguishing between
different orbifold types of pending vertices), we
associate a {\em mapping class group} operation to this morphism
therefore passing from the groupoid of morphisms to the group of
modular transformations.

\subsubsection{Changing the decoration (spiraling direction)}

The last class of mapping class group transformations pertains to changing the ``decoration,'' i.e.,
the sign of the hole perimeter:
\be
\label{loopinvert}
{\psset{unit=0.5}
\begin{pspicture}(-5,-3)(7,1)
\pcline[linewidth=2pt,linecolor=blue](-6,-0.5)(-4,-0.5)
\pcline[linewidth=2pt,linecolor=blue](-6,-1.5)(-4,-1.5)
\psbezier[linewidth=2pt,linecolor=blue](-4,-0.5)(-3,1)(-1,1)(-1,-1)
\psbezier[linewidth=2pt,linecolor=blue](-4,-1.5)(-3,-3)(-1,-3)(-1,-1)
\psbezier[linewidth=2pt,linecolor=blue](-3.2,-1)(-2.4,-0.2)(-2,-0.3)(-2,-1)
\psbezier[linewidth=2pt,linecolor=blue](-3.2,-1)(-2.4,-1.8)(-2,-1.7)(-2,-1)
\rput(-5.5,0.2){\makebox(0,0){$Y$}}
\rput(-1.5,1){\makebox(0,0){$P$}}
\pcline[linewidth=2pt]{<->}(0,-1)(2,-1)
\pcline[linewidth=2pt,linecolor=blue](3,-0.5)(5,-0.5)
\pcline[linewidth=2pt,linecolor=blue](3,-1.5)(5,-1.5)
\psbezier[linewidth=2pt,linecolor=blue](5,-0.5)(6,1)(8,1)(8,-1)
\psbezier[linewidth=2pt,linecolor=blue](5,-1.5)(6,-3)(8,-3)(8,-1)
\psbezier[linewidth=2pt,linecolor=blue](5.8,-1)(6.6,-0.2)(7,-0.3)(7,-1)
\psbezier[linewidth=2pt,linecolor=blue](5.8,-1)(6.6,-1.8)(7,-1.7)(7,-1)
\rput(3.5,0.2){\makebox(0,0){$Y+P$}}
\rput(7.5,1){\makebox(0,0){$-P$}}
\rput(10,0){\makebox(0,0){.}}
\pcline[linecolor=red, linewidth=1pt](-6,-.7)(-3.9,-.7)
\pcline[linecolor=red, linewidth=1pt](-6,-1.3)(-3.9,-1.3)
\psbezier[linecolor=red, linewidth=1pt]{->}(-3.9,-.7)(-3,.7)(-1.3,.8)(-1.3,-1)
\psbezier[linecolor=red, linewidth=1pt](-3.9,-1.3)(-3,-2.7)(-1.3,-2.8)(-1.3,-1)
\pcline[linecolor=red, linewidth=1pt](3,-.7)(5.1,-.7)
\pcline[linecolor=red, linewidth=1pt](3,-1.3)(5.1,-1.3)
\psbezier[linecolor=red, linewidth=1pt]{->}(5.1,-.7)(6,.7)(7.7,.8)(7.7,-1)
\psbezier[linecolor=red, linewidth=1pt](5.1,-1.3)(6,-2.7)(7.7,-2.8)(7.7,-1)
\end{pspicture}
}
\ee
That this transformation preserves the Poisson brackets is obvious because the variable $P$ Poisson commutes with all
other variables, whereas the preservation of geodesic functions follows from two matrix equalities:
\bea
&&X_YLX_PLX_Y=X_{Y+P}LX_{-P}LX_{Y+P},
\nonumber
\\
&&X_YRX_PRX_Y=X_{Y+P}RX_{-P}RX_{Y+P}.
\nonumber
\eea

We can therefore extend the mapping class group of ${\mathfrak T}^H_{g,s,r}$ by adding
symmetries between sheets of the $2^s$-ramified covering of the ``genuine'' (nondecorated)
Teichm\"uller space ${\mathfrak T}_{g,s,r}$.

We can summarize as follows.

\begin{theorem}
The whole mapping class group of $\Sigma_{g,s,r}$ is generated by morphisms described by
Lemmas~\ref{lem-abc} and~\ref{lem-pending1} and formula~(\ref{loopinvert}).
\end{theorem}

\subsection{Quantum MCG transformations}
Recall that a quantization of a Poisson manifold equivariant w.r.t. a
discrete group action is a family of $*$-algebras ${\cal A}^\hbar$ depending
smoothly on a positive real parameter $\hbar$, acting on~$G$ by outer
automorphisms and
having the following relation to the Poisson manifold.

{\bf 1.} For $\hbar=0$, the algebra is isomorphic as a $G$-module to the
$*$-algebra of complex-valued function on the Poisson manifold.

{\bf 2.} The Poisson bracket on ${\cal A}^0$ given by $\{a_1, a_2\} =
\lim_{\hbar \rightarrow 0}\frac{[a_1,a_2]}{\hbar}$ coincides with the one
generated by the Poisson structure of the manifold.

We now quantize a Teichm\"uller space
${\mathfrak T}^{H}_{g,s,r}$ equivariantly w.r.t.\ the mapping class group action.

Let ${\mathfrak T}^\hbar(\Gamma_{g,s,r})$ be a
$*$-algebra generated by the generator $Z_\alpha^\hbar$
(one generator per one unoriented edge $\alpha$)
and relations
\be
[Z^\hbar_\alpha, Z^\hbar_\beta ] = 2\pi i\hbar\{Z_\alpha, Z_\beta\}
\label{qq}
\ee
with the $*$-structure
\be
(Z^\hbar_\alpha)^*=Z^\hbar_\alpha.
\ee
Here $z_\alpha$  and $\{\cdot,\cdot\}$ stand for the respective
coordinate functions on
the classical Teichm\"uller space and the Weil--Petersson Poisson bracket on
it. Note that according to formula (\ref{WP-PB}), the
right-hand side of (\ref{qq})
is merely a constant which may take only five values $0$, $\pm
2\pi i \hbar$, $\pm 4 \pi i \hbar$.

For the notation simplicity in what follows we omit the superscript $\hbar$ for the quantum operators;
the classical or quantum nature of the object will be always clear from the context.

\subsubsection{Quantum flip morphisms for inner edges}

It was proved in \cite{ChF1} that the {\em quantum flip morphisms}
\bea
\{A,B,C,D,Z\}&\to&\{A+\phi^\hbar(Z),B-\phi^\hbar(-Z),C+\phi^\hbar(Z),
D-\phi^\hbar(-Z),-Z\}\nonumber\\
&:=&\{\tilde A,\tilde B,\tilde C,\tilde D,\tilde Z\},
\label{q-mor}
\eea
where $A$, $B$, $C$, $D$, and~$Z$ are as in Fig.~\ref{fi:flip} and
$\phi^\hbar(x)$ is the real function of one real variable,
\begin{equation} \label{phi}
\phi^\hbar(z) =
-\frac{\pi\hbar}{2}\int_{\Omega} \frac{e^{-ipz}}{\sinh(\pi p)\sinh(\pi \hbar
p)}dp,
\end{equation}
(the contour $\Omega$ goes along the real axis bypassing the
singularity at the origin from above) satisfy the standard two-, four-, and five-term relations.
The quantum dilogarithm function $\phi^\hbar(z)$ was introduced in this context by Faddeev in
\cite{Faddeev} and used in \cite{FK} for constructing quantum MCG transformations for the Liouville model.

As in the classical case, the center of the algebra is generated by the sums
$Z^\hbar_f :=\sum_{\alpha \in f}{Z^\hbar_\alpha}$ ranging all edges $\alpha$
surrounding a given face $f$.

Important properties of the quantum flip morphisms are
\begin{itemize}
\item In the limit $\hbar \rightarrow 0$, morphism
{\rm(\ref{q-mor})} coincides with the classical morphism, that is,
$$\lim_{\hbar \rightarrow 0}\phi^\hbar(z) = \log(e^z + 1).$$
\item The morphism of the commutation relations is ensured by the relation
$$
\phi^\hbar(z)-\phi^\hbar(-z)=z.
$$
\item The Hermiticity of the transformed operators is ensured by the conjugation condition
$$
\overline{\phi^\hbar(z)} = \phi^\hbar(\overline{z}).
$$
\item Eventually, we have the double quasi-periodicity conditions responsible for the satisfaction of
the five-term relation:
\bea
&{}&
\phi^\hbar(z+i\pi\hbar)-\phi^\hbar(z-i\pi\hbar) = \frac{2\pi i
\hbar}{e^{-z}+1},\label{phi-period1}\\
&{}&
\phi^\hbar(z+i\pi)-\phi^\hbar(z-i\pi) = \frac{2\pi
i}{e^{-z/\hbar}+1}.\label{phi-period2}
\eea
\item The function $\phi^\hbar(z)$ is meromorphic with poles at the
points $\{\pi i(m+ n\hbar)|m,n \in {\bf N}\}$
and $\{-\pi i (m+n\hbar)|m,n \in {\bf N}\}$.
\end{itemize}

We also have conditions of the modular double, $\frac{1}{\hbar}\phi^\hbar(z) = \phi^{1/\hbar}(z/\hbar)$ and
the first quasi-periodicity condition (\ref{phi-period1}), which are irrelevant when we discuss geodesic functions
on which half of the modular double generators act trivially.

\begin{remark}
Note that exponentiated
algebraic elements $U_i = e^{\pm Z_i}$,
which obey homogeneous commutation relations
$q^nU_iU_j=U_jU_i q^{-n}$ with $[X_i,X_j]=2ni\pi  \hbar $ and $q:=e^{i\pi\hbar}$
transform as rational functions: for example,
\bea
&{}&
e^{A+\phi^\hbar(Z)}=\exp\left(\frac{1}{2\pi i
\hbar}\int_{A}^{A+2i\pi \hbar}\phi^\hbar(z)dz\right)e^{A}=\nonumber\\
&{}&\qquad\qquad
=\exp\left(\frac{1}{2i\pi \hbar}\int_{-\infty}^{A}(\phi^\hbar(z+2i\pi \hbar)
-\phi^\hbar(z))dz\right) e^{A}=\nonumber\\
&{}&\qquad\qquad
= \exp\left(\int_{-\infty}^{A}\frac{dz}{e^{-z-\pi i \hbar} +1}\right)e^A =(1+qe^Z)e^A,
\label{trans}
\eea
where we have used the standard formula
$$
e^{A+\Phi(B)} =
\exp\left\{\frac{1}{[A,B]}\int_B^{B+[A,B]}\Phi(z)dz\right\}e^A,
$$
which is valid for all~$A$ and~$B$ such that the commutator $[A,B]$ is a
nonzero scalar.
\end{remark}

\subsubsection{Quantum flip morphisms for pending edges}

We now present the new formula, which is the quantization of the MCG transformation in
Fig.~\ref{fi:interchange-p-dual}.

\begin{lemma} \label{lem-pending-quantum}
The transformation~in Fig.~\ref{fi:interchange-p-dual}
\be
\label{morphism-pending-quantum}
\{\tilde A,\tilde B,\tilde Z\}:= \{A+\phi^\hbar(Z+i\varphi)+\phi^\hbar(Z-i\varphi), B-\phi^\hbar(-Z+i\varphi)-\phi^\hbar(-Z-i\varphi),-Z\},
\ee
with $\phi^\hbar(x)$ from (\ref{phi}) and $w=e^{i\varphi}+e^{-i\varphi}$,
is the morphism of the quantum $*$-algebra
${\mathfrak T}^\hbar_{g,s,r}$.
\end{lemma}

Verifying the preservation of the commutation relations is simple, the most difficult part is to verify
the preservation of geodesic function operators. We devote a special section to follow to the verification
of this condition.

\section{Operator ordering for quantum geodesics}\label{s:ordering}

\subsection{Transformation properties of quantum geodesic functions}\label{ss:QGF}

The quantum analogues of {\em matrix} relations (\ref{mutation1})--(\ref{mutation3}) were found in
\cite{ChP1}. Here and hereafter, for the rest of the paper, we assume that the ordering of
quantum operators in a product is {\em natural}, i.e.,
it is determined by the order of matrix multiplication itself. Then,
amazingly enough, {\em all four} entries of the corresponding $2\times 2$-matrices transform uniformly:
applying the quantum MCG transformation to curves 1,2, and 3 in Fig.~\ref{fi:flip}, we obtain
the respective quantum matrix relations:
\bea
X_D RX_ZRX_A &=&q^{1/4}X_{\tilde D}RX_{\tilde A},\label{curve1}\\
X_D RX_Z LX_B&=& X_{\tilde D} RX_{\tilde Z} LX_{\tilde B},\label{curve2}\\
X_D L X_C &=&q^{1/4}X_{\tilde D}L X_{\tilde Z} L X_{\tilde C}.\label{curve3}
\eea
We prove here the third relation; the remaining relations can be verified analogously.
In the left-hand side, we obtain, after reducing the matrix elements to the Weyl-ordered form, that
\bea
X_{\tilde D}L X_{\tilde Z} L X_{\tilde C}
&=&\left[
\begin{array}{cc}
e^{\tilde D/2-\tilde C/2-\tilde Z/2}+e^{\tilde D/2-\tilde C/2+\tilde Z/2} &
-q^{-1/2}e^{\tilde D/2+\tilde C/2+\tilde Z/2} \\
-q^{-1/2}e^{-\tilde D/2-\tilde C/2-\tilde Z/2} & 0 \\
\end{array}
\right]
\nonumber\\
&=&\left[\begin{array}{cc}
e^{D/2- C/2} & -q^{-1/2}e^{D/2+C/2} \\
-q^{-1/2}e^{- D/2-C/2} & 0 \\
\end{array}
\right]\nonumber,
\eea
where we have used the relation $\tilde D+\tilde C+\tilde Z=D+C$ and the chain of equalities
\bea
&{}&e^{\tilde D/2-\tilde C/2-\tilde Z/2}+e^{\tilde D/2-\tilde C/2+\tilde Z/2}
=e^{Z+D/2-C/2-\phi^\hbar(Z)}+e^{D/2-C/2-\phi^\hbar(Z)}\nonumber\\
&{}&=\exp\Bigl\{\frac{1}{2\pi i\hbar}\int_{Z}^{Z-2\pi i \hbar}[-\xi+\phi^\hbar(\xi)]d\xi
\Bigr\}e^{D/2-C/2}+\exp\Bigl\{\frac{1}{2\pi i\hbar}\int_{Z}^{Z-2\pi i \hbar}\phi^\hbar(\xi)d\xi
\Bigr\}e^{D/2-C/2}\nonumber\\
&{}&=\left[e^{Z-i\pi \hbar}+1\right]
\exp\Bigl\{\frac{1}{2\pi i\hbar}\int_{-\infty}^{Z}\frac{-2\pi i\hbar d\xi}{e^{-\xi+i\pi \hbar}+1}
\Bigr\}e^{D/2-C/2}\nonumber\\
&{}&=\left[e^{Z-i\pi \hbar}+1\right]\frac{1}{e^{Z-i\pi \hbar}+1}e^{D/2-C/2}=e^{D/2-C/2}.\nonumber
\eea
In the right-hand side, we have
$$
X_D L X_C =\left[
             \begin{array}{cc}
               q^{1/4}e^{D/2-C/2} & -q^{-1/4}e^{C/2+D/2} \\
               -q^{-1/4}e^{-C/2-D/2} & 0 \\
             \end{array}
           \right]
$$
and, comparing powers of $q$ standing with each of the four matrix entries, we come to the
matrix equality (\ref{curve3}).

We now formulate the new lemma concerning quantum versions of MCG transformations
(\ref{pend-mutation1})--(\ref{pend-mutation3}).

\begin{lemma}\label{lem-quantum-geod}
We have the following quantum matrix relations:
\bea
X_{A}LX_Z(F_p\Omega) X_ZLX_{B}&=&q^{-1}
X_{{\tilde A}}RX_{\tilde Z}(-\Omega) X_{\tilde Z}RX_{{\tilde B}},
\label{pend-quantum-mutation1}\\
X_{A}LX_Z\Omega X_ZRX_{A}&=&
X_{{\tilde A}}RX_{\tilde Z}(-\Omega)X_{\tilde Z}LX_{{\tilde A}},
\label{pend-quantum-mutation2}\\
X_{B}RX_Z\Omega X_ZLX_{B}&=&
X_{{\tilde B}}LX_{\tilde Z}(-\Omega)X_{\tilde Z}RX_{{\tilde B}}.
\label{pend-quantum-mutation3}
\eea
\end{lemma}

The {\em proof} is again a direct calculation.

\begin{remark}\label{rem-quantum}
A closer look on transformation laws (\ref{curve1})--(\ref{curve3}) and
(\ref{pend-quantum-mutation1})--(\ref{pend-quantum-mutation3}) yields a ``mnemonic'' rule for the
powers of $q$: the ``proper'' power of $q$ is always one fourth of the difference between the number of
matrices $R$ and $L$ in the matrix chain. Although for the general geodesic functions this property presumably does
not lead to meaningful answers coming into contradiction with the Hermiticity condition,
for a special class of geodesic functions in the next section it enables us to solve
the problem of quantum ordering completely.
\end{remark}

\subsection{Quantum ordering for special class of geodesic functions}\label{ss:invariance}

Quantum matrix MCG relations turn out to be very useful when we consider {\em special} geodesic functions that correspond to geodesics going around exactly two orbifold points or holes. Fortunately,
many important geodesic functions fall into
this class, including those appearing in the description of Poisson leaves of twisted Yangian systems of $A_n$ and
$D_n$ types \cite{ChM} and those pertaining to the description of Painlev\'e VI equation \cite{ChM-D4}. Our aim in this and the next section is to introduce a proper quantum ordering not only for quantum geodesic functions but for
more general monodromy matrices of Schlesinger systems (we then obtain geodesic functions as traces of products of
these monodromy matrices).

We first take one of these points (the zeroth orbifold point) as a root, let $S:=Z_0$
be the variable of the pending edge terminating at this point, and consider
$2\times 2$-{\em monodromy matrices} $M_i$ corresponding to going around an $i$th orbifold point. These matrices
have the structure
\be
M_i=X_SRX_{Y_1}L\cdots R X_{Y_p}L X_{Z_i}F_{\omega_i}X_{Z_i}R X_{Y_p}L\cdots R X_{Y_1}LX_S, \ i\ne0,\quad
M_0=F_{\omega_0},
\label{monodromy}
\ee
where $Z_i$ is the variable of the pending edge terminating at the $i$th orbifold point and $Y_1,\dots,Y_p$ are
variables of the intermediate (not necessarily internal) edges.
In the general situation we require the path joining the root and the
orbifold point neither to be unique (this is true only for fat graphs without loops) nor to have no self-intersections. Whatever is this path, from the transformation properties of quantum geodesic functions
we can immediately conclude that any product of quantum matrices of the form (\ref{monodromy}) is invariant under
all the quantum MCG transformations that do not mutate the root edge. For a path without self-intersections,
using quantum MCG transformations
we can reduce any string of quantum matrices (\ref{monodromy}) to
\be
M_i=X_{\tilde S}LX_{\tilde Z_i}F_{\omega_i}X_{\tilde Z_i}RX_{\tilde S}
=\left(
\begin{array}{cc}
qe^{-\tilde Z_i}+\omega_i &
-e^{-\tilde Z_i+\tilde S}-e^{\tilde Z_i+\tilde S}-\omega_ie^{\tilde S} \\
e^{-\tilde Z_i-\tilde S} & -q^{-1}e^{-\tilde Z_i} \\
\end{array}
\right),\quad [\tilde Z_i,\tilde S]=2\pi i\hbar.
\label{monodromy2}
\ee
Introducing the notation
\be
a_i:=e^{-\tilde Z_i},\quad b_i=e^{-\tilde Z_i+\tilde S}+e^{\tilde Z_i+\tilde S}+\omega_ie^{\tilde S},
\quad c_i:=e^{-\tilde Z_i-\tilde S},
\label{abc-i}
\ee
where $a_i$, $b_i$, and $c_i$ are Hermitian operators,
we eventually construct from the entries of the monodromy matrix (\ref{monodromy2}) (or (\ref{monodromy}))
the quantum geodesic function that is invariant also w.r.t. mutations of the root edge (at the end of which we
set the orbifold point with the parameter $\omega_0$):
\be
G_{0,i}=b_i+c_i+\omega_0 a_i.
\label{G0i}
\ee
We have thus completely solved the problem of quantum ordering for the product of quantum
edge matrices entering the geodesic functions corresponding to geodesics separating exactly two orbifold points
from the rest of the Riemann surface.

\section{Quantum relations for monodromy matrices}\label{s:quantum-monodromy}

We know that for some special cases of orbifold Riemann surfaces,
the algebra of geodesic functions can be closed. We demonstrate that, in the very same cases, we are also able
to close the algebras of elements the Fuchsian group, or monodromy matrices.
In this paper, we consider only two such cases: the first one is the $A_n$-algebra of a sphere
with one hole and with $n$ $\mathbb Z_2$-orbifold points (in this case, all $\omega_i$ are zeros,
and $X_Z FX_Z=X_{2Z}$); the second one is the
algebra of geodesic functions on a sphere with four holes/orbifold points of arbitrary orders. We leave
to a separate publication the most
interesting but also the most difficult case of $D_n$-algebras related to reflection equations with the spectral parameter.

\subsection{$A_n$-algebra of monodromy matrices}\label{ss:An-monodromy}

\begin{figure}[tb]
\begin{pspicture}(-6,-1)(4,2)
\rput{90}(4,0){
\parametricplot[linecolor=red,linewidth=1pt]{0}{8}{1 2 t mul add -1 exp 1 add 60 t mul sin mul 1 2 t mul add -1 exp 1 add 60 t mul cos mul }}
\rput{270}(0,0){\parametricplot[linecolor=red,linewidth=1pt]{0}{3}{1 2 t mul add -1 exp 1 add 60 t mul sin mul 1 2 t mul add -1 exp 1 add 60 t mul cos mul }}
\rput{90}(4,0){
\parametricplot[linecolor=red,linewidth=1pt]{0}{8}{1 2 t mul add -1 exp 1.05 add 49 t exp 2 mul -98 t mul add 154 add -1 exp 484 mul add 60 t mul sin mul 1 2 t mul add -1 exp 1.05 add  49 t exp 2 mul -98 t mul add 154 add -1 exp 484 mul add  60 t mul cos mul }}
\psbezier[linecolor=red,linewidth=1pt,linestyle=dashed](0,0)(-0.2,1.2)(-0.2,1.5)(0,2.7)
\psbezier[linecolor=red,linewidth=1pt,linestyle=dashed](0,0)(1,1.6)(1,1.6)(2,2.7)
\psbezier[linecolor=red,linewidth=1pt,linestyle=dashed](0,0)(1.5,1.2)(1.5,1.2)(3,2)
\rput(0.6,1.7){\makebox(0,0){{\tcr{\tiny$\bullet$}}}}
\rput(0.4,1.7){\makebox(0,0){{\tcr{\tiny$\bullet$}}}}
\rput(0.2,1.7){\makebox(0,0){{\tcr{\tiny$\bullet$}}}}
\rput(-.2,-.4){\makebox(0,0){$p_0$}}
\rput(0,0){\makebox(0,0){\tcb{$\bullet$}}}
\rput(0.1,3.1){\makebox(0,0){$p_{n-1}$}}
\rput(0,2.7){\makebox(0,0){\tcb{$\bullet$}}}
\rput(0.8,2.7){\makebox(0,0){$\cdot$}}
\rput(1,2.7){\makebox(0,0){$\cdot$}}
\rput(1.2,2.7){\makebox(0,0){$\cdot$}}
\rput(2.1,3.1){\makebox(0,0){$p_2$}}
\rput(2,2.7){\makebox(0,0){\tcb{$\bullet$}}}
\rput(3.2,2.4){\makebox(0,0){$p_1$}}
\rput(3,2){\makebox(0,0){\tcb{$\bullet$}}}
\end{pspicture}
\caption{\small
The domain bounded by an infinite geodesic line (a side of a triangle in an ideal triangular decomposition of the Riemann surface) inside which all geodesics $\gamma_{0,i}$ joining orbifold points $p_0$ and $p_i$, $i=1,\dots,n-1$,
and constituting the $A_n$ system must be contained.
}
\label{fi:domain}
\end{figure}
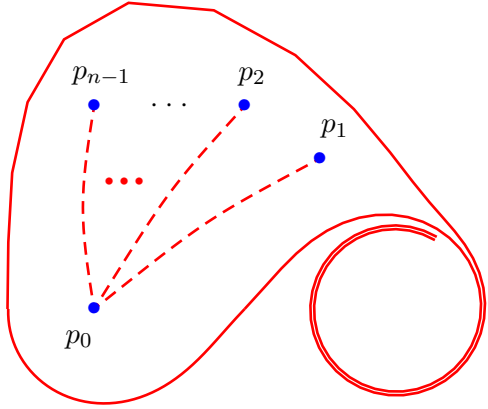

From the geometry side, we consider a pattern in which $n$ ${\mathbb Z}_2$-orbifold points are separated
from the rest of the Riemann surface by an infinite geodesic line (a side of an ideal triangle from an ideal
triangle decomposition of a finite part of the Riemann surface) whose both ends wind to the perimeter line of
the same hole as shown in Fig.~\ref{fi:domain}. The relevant geodesic lines joining the root orbifold point $p_0$
with other $n-1$ orbifold points $p_i$, $i=1,\dots,n-1$, and constituting the $A_n$ set must then have no (self)intersections and must be contained inside this (simply connected) domain. We then have a natural linear
ordering of these $n-1$ orbifold points ordered counterclockwise.

Given two quantum monodromy matrices $M_1$ and $M_2$ corresponding to paths joining $p_0$ with
two other $\mathbb Z_2$-orbifold points, one of which ($p_1$) precedes the other ($p_2$) in the
corresponding linear ordering,
we can always, using the quantum MCG
transformations from the previous sections, reduce them to the case in which (see Fig.~\ref{fi:M1-M2})
\bea
M_1&=&X_S LX_X L X_{2Z} R X_X R X_S,\label{M1}\\
M_2&=&X_S LX_X R X_{2Y} L X_X R X_S,\label{M2}
\eea
with the nontrivial quantum commutation relations
\bea
&{}&e^X e^S=qe^{X+S},\quad e^Xe^Y=qe^{X+Y},\nonumber\\
&{}&e^Y e^Z=qe^{Y+Z},\quad e^Ze^X=qe^{Z+X}.\nonumber
\eea
The both matrices $M_i$ have the same structure,
\be
M_i=\left[
      \begin{array}{cc}
        qa_i & -b_i \\
        c_i & -q^{-1}a_i \\
      \end{array}
    \right];
    \label{ABC-new}
\ee
in which we have explicit formulas for $a_i$, $b_i$, and $c_i$:
\bea
a_1&=&e^{-X-Z}+e^{-Z}\nonumber\\
b_1&=&e^{-X-Z+S}+(q+q^{-1})e^{-Z+S}+e^{X-Z+S}+e^{X+Z+S}\nonumber\\
c_1&=&e^{-X-Z-S},\nonumber
\eea
and
\bea
a_2&=&e^{-X-Y}+e^{-X+Y}+e^{Y}\nonumber\\
b_2&=&e^{-X-Y+S}+(q+q^{-1})e^{Y+S}+e^{-X+Y+S}+e^{X+Y+S}\nonumber\\
c_2&=&e^{-X-Y-S}+e^{-X+Y-S}.\nonumber
\eea
In what follows, all $a_i$, $b_i$, and $c_i$ are Hermitian operators.

\begin{figure}[tb]
\begin{pspicture}(-4,-1.5)(4,1.5){
\newcommand{\MONE}{%
{\psset{unit=0.7}
\psline[linewidth=18pt,linecolor=blue](-2,0)(0,0)
\psline[linewidth=18pt,linecolor=blue](0,0)(1,1.5)
\psline[linewidth=18pt,linecolor=blue](0,0)(1,-1.5)
\psline[linewidth=18pt,linecolor=blue](-3,-1.5)(-2,0)
\psline[linewidth=18pt,linecolor=blue](-2,0)(-2.5,.75)
\psline[linewidth=14pt,linecolor=white](-2,0)(0,0)
\psline[linewidth=14pt,linecolor=white](0,0)(1,1.5)
\psline[linewidth=14pt,linecolor=white](0,0)(1,-1.5)
\psline[linewidth=14pt,linecolor=white](-3,-1.5)(-2,0)
\psline[linewidth=14pt,linecolor=white](-2,0)(-2.5,.75)
%
\rput(1,-1.5){\pscircle[linewidth=2pt,linecolor=red]{0.25}}
\rput(0,0){\pscircle[linewidth=2pt,linecolor=red]{0.25}}
\rput(-2,0){\pscircle[linewidth=2pt,linecolor=red]{0.25}}
\psline[linewidth=10pt,linecolor=red](-2,0)(0,0)
\psline[linewidth=10pt,linecolor=red](0,0)(1,-1.5)
\psline[linewidth=10pt,linecolor=red](-3,-1.5)(-2,0)
\psline[linewidth=6pt,linecolor=white](-2,0)(0,0)
\psline[linewidth=6pt,linecolor=white](0,0)(1,-1.5)
\psline[linewidth=6pt,linecolor=white](-3,-1.5)(-2,0)
\rput(1,1.5){\pscircle*{0.05}}
\rput(1,-1.5){\pscircle*{0.05}}
\rput(0.4,1.5){\makebox(0,0)[rb]{$Y$}}
\rput(-1,0.7){\makebox(0,0)[cb]{$X$}}
\rput(0.4,-1.5){\makebox(0,0)[rt]{$Z$}}
\rput(-2.4,-1.5){\makebox(0,0)[lt]{$S$}}
\rput(-3,-1.5){\makebox(0,0)[rt]{$M_1$}}
}
}
\newcommand{\MTWO}{%
{\psset{unit=0.7}
\psline[linewidth=18pt,linecolor=blue](-2,0)(0,0)
\psline[linewidth=18pt,linecolor=blue](0,0)(1,1.5)
\psline[linewidth=18pt,linecolor=blue](0,0)(1,-1.5)
\psline[linewidth=18pt,linecolor=blue](-3,-1.5)(-2,0)
\psline[linewidth=18pt,linecolor=blue](-2,0)(-2.5,.75)
\psline[linewidth=14pt,linecolor=white](-2,0)(0,0)
\psline[linewidth=14pt,linecolor=white](0,0)(1,1.5)
\psline[linewidth=14pt,linecolor=white](0,0)(1,-1.5)
\psline[linewidth=14pt,linecolor=white](-3,-1.5)(-2,0)
\psline[linewidth=14pt,linecolor=white](-2,0)(-2.5,.75)
%
\rput(1,1.5){\pscircle[linewidth=2pt,linecolor=red]{0.25}}
\rput(0,0){\pscircle[linewidth=2pt,linecolor=red]{0.25}}
\rput(-2,0){\pscircle[linewidth=2pt,linecolor=red]{0.25}}
\psline[linewidth=10pt,linecolor=red](-2,0)(0,0)
\psline[linewidth=10pt,linecolor=red](0,0)(1,1.5)
\psline[linewidth=10pt,linecolor=red](-3,-1.5)(-2,0)
\psline[linewidth=6pt,linecolor=white](-2,0)(0,0)
\psline[linewidth=6pt,linecolor=white](0,0)(1,1.5)
\psline[linewidth=6pt,linecolor=white](-3,-1.5)(-2,0)
\rput(1,1.5){\pscircle*{0.05}}
\rput(1,-1.5){\pscircle*{0.05}}
\rput(0.4,1.5){\makebox(0,0)[rb]{$Y$}}
\rput(-1,0.7){\makebox(0,0)[cb]{$X$}}
\rput(0.4,-1.5){\makebox(0,0)[rt]{$Z$}}
\rput(-2.4,-1.5){\makebox(0,0)[lt]{$S$}}
\rput(-3,-1.5){\makebox(0,0)[rt]{$M_2$}}
}
}
\rput(2.5,0){\MONE}
\rput(7.5,0){\MTWO}
}
\end{pspicture}
\caption{\small
The graphical representation for the monodromy matrices $M_1$ and $M_2$.
}
\label{fi:M1-M2}
\end{figure}

Note that the matrices $\tilde M_i:= F_0 M_i$ corresponding to closed geodesics ($F_0$ is the matrix with $\omega=0$)
satisfy the $U_q(\mathfrak{sl}_2)$ relations. We have
\be
\begin{array}{ll}
qa_ib_i-q^{-1}b_ia_i=0&\\
q^{-1}a_ic_i-qc_ia_i=0,&\qquad i=1,2 \\
b_ic_i-c_ib_i=(q^2-q^{-2})a_i^2\quad\hbox{in fact}\quad b_ic_i=q^2a_i^2+1,&\\
\end{array}
\label{sl2q}
\ee
where 1 is the unit element. So, we are in the situation of a quantum algebra with the unity.
The above matrices satisfy the quantum relation
\be
M_i^2=-{\mathbb E}.
\label{M-square}
\ee

Interestingly, we can find the {\em complete set} of quantum algebraic relations between entries of $M_1$ and $M_2$, thus leading to a nice algebraic structure
on $U_q(\mathfrak{sl}_2)\times U_q(\mathfrak{sl}_2)$ (we present it in the form closest to the
$R$-matrix one):
\bea
&{}& q^{-1}b_1b_2=qb_2b_1,\nonumber\\
&{}& q^{-1}c_1c_2=qc_2c_1,\nonumber\\
&{}& a_1b_2=b_2a_1,\nonumber\\
&{}& b_1a_2=a_2b_1+(q^2-q^{-2})a_1b_2=a_2b_1+(q^2-q^{-2})b_2a_1,\nonumber\\
&{}& c_1a_2=a_2c_1,\label{algebraM1-M2-r-matrix}\\
&{}& a_1c_2 =c_2a_1+(q^2-q^{-2})c_1a_2=c_2a_1+(q^2-q^{-2})a_2c_1,\nonumber\\
&{}& qc_1b_2=q^{-1}b_2c_1,\nonumber\\
&{}&a_1a_2=a_2a_1+(1-q^{-2})b_2c_1=a_2a_1+(q^2-1)c_1b_2,\nonumber\\
&{}&qb_1c_2+q[q^{-2}-q^{2}]a_1a_2=q^{-1}c_2b_1-q^{-1}[q^{-2}-q^{2}]a_2a_1.\nonumber
\eea

Since we are in a general position case, we can formulate the algebra for {\em any} system of
$2\times 2$-matrices $M_i$, $i=1,\dots,n-1$
such that the paths between the root and all orbifold points are drawn on a
tree-like fat subgraph; for simplicity we assume that there exists such a choice of
a fat graph in which all these paths are graph-simple, i.e., in representation (\ref{monodromy}) $Y_q\ne Y_r$ for
$1\le q<r\le p$. We have then the natural linear ordering of the orbifold points from the root of the tree counterclockwise. The algebra of entries $a_i$, $b_i$, $c_i$ is given by (\ref{sl2q})
For $i<j$, we then have the algebra (\ref{algebraM1-M2-r-matrix}) in which we merely substitute
$i$ and $j$ for the respective indices 1 and 2 of the elements of monodromy matrices:
\bea\label{sl2qxsl2q}
&{}& q^{-1}b_ib_j=qb_jb_i,\nonumber\\
&{}& q^{-1}c_ic_j=qc_jc_i,\nonumber\\
&{}& a_ib_j=b_ja_i,\nonumber\\
&{}& b_ia_j=a_jb_i+(q^2-q^{-2})a_ib_j=a_jb_i+(q^2-q^{-2})b_ja_i,\nonumber\\
&{}& c_ia_j=a_jc_i,\label{algebraMi-Mj}\\
&{}& a_ic_j =c_ja_i+(q^2-q^{-2})c_ia_j=c_ja_i+(q^2-q^{-2})a_jc_i,\nonumber\\
&{}& qc_ib_j=q^{-1}b_jc_i,\nonumber\\
&{}&a_ia_j=a_ja_i+(1-q^{-2})b_jc_i=a_ja_i+(q^2-1)c_ib_j,\nonumber\\
&{}&qb_ic_j+q[q^{-2}-q^{2}]a_ia_j=q^{-1}c_jb_i-q^{-1}[q^{-2}-q^{2}]a_ja_i.\nonumber
\eea

\begin{prop} The abstract quantum algebra for $a_i,b_i,c_i$ , $i=1,\dots,n-1$ defined by the relations (\ref{sl2q}) and (\ref{sl2qxsl2q}) satisfies the quantum Jacobi relation. \end{prop}

{\em Proof.} This is a consequence of the fact that the formulae defining  $a_i,b_,c_i$  , $i=1,\dots,n-1$ in terms of the Darboux coordinates are invertible in a large open set.

\vskip 2mm

It is interesting to compare these algebras with the known algebras of traces of products of monodromy
matrices like in \cite{ChM}. Recalling that the monodromy matrix around the root is $M_0:=F_0=\left[
                                                           \begin{array}{cc}
                                                             0 & 1 \\
                                                             -1 & 0 \\
                                                           \end{array}
                                                         \right]$,
the obvious candidates for $G_{0,i}=\tr M_iM_0$ are linear combinations
\be
G_{0,i}=b_i+c_i.
\label{G-lin}
\ee
A more subtle question is what is $G_{i,j}$, which is the quantum analogue of the trace of product of $M_i$ and $M_j$? Here, we encounter the problem of quantum ordering. From the comparison with the known answer ($G_{1,2}=e^{Y+Z}+e^{-Y+Z}+e^{-Y-Z}$) after some algebra we obtain
\be
G_{i,j}=qb_ic_j+q^3 c_ib_j-(q^3+q)a_ia_j, \quad 1\le i<j\le n-1.
\label{Gij}
\ee

Let us verify the Hermiticity condition:
\bea
G^*_{1,2}&=&q^{-1}c_2b_1+q^{-3}b_2 c_1-(q^{-3}+q^{-1})a_2a_1\nonumber\\
&=&qb_1c_2-[q^2-q^{-2}][qa_1a_2+q^{-1}a_2a_1] +q^{-1}c_1b_2-(q^{-3}+q^{-1})a_2a_1\nonumber\\
&=&qb_1c_2 +q^{-1}c_1b_2-q^3a_1a_2-qa_2a_1+q^{-1}[a_1,a_2]\nonumber\\
&=&qb_1c_2 +q^{-1}c_1b_2-q^3a_1a_2+q^{-1}[q-q^{-1}]qc_1b_2-qa_2a_1\nonumber\\
&=&qb_1c_2 +q c_1b_2-q^3a_1a_2-qa_2a_1\nonumber\\
&=&qb_1c_2 +q c_1b_2-q^3a_1a_2-qa_1a_2+q^2[q-q^{-1}]c_1b_2=G_{1,2}.\nonumber
\eea
Using Matematica NCAlgebra package \cite{NCA} we have verified all the relations of the Nelson--Regge algebra~\cite{NR},~\cite{NRZ},~\cite{Ugaglia}:
\bea
&{}&[G_{i,j},G_{k,l}]=0\quad\hbox{for}\quad 0\le i<j<k<l\le n-1\ \hbox{and}\ 0\le i<k<l<j\le n-1;\nonumber\\
&{}&[G_{i,k},G_{j,l}]=[q^2-q^{-2}][G_{i,j}G_{k,l}- G_{i,l}G_{j,k}]\quad\hbox{for}
\quad 0\le i<j<k<l\le n-1;\nonumber\\
&{}&q G_{i,j}G_{j,k}-q^{-1}G_{j,k}G_{i,j}=[q^2-q^{-2}]G_{i,k}\quad\hbox{for}\quad 0\le i<j<k\le n-1,\ \hbox{etc}.\nonumber
\eea

\begin{remark}\label{rem:parameter-counting}
Note that each monodromy matrix $M_i$, $i=1,\dots,n-1$,
brings two independent quantum operators because we always have
one relation (the last formula in (\ref{sl2q})) on three additional operators. The same is true for the
geometric system: our rooted fat subgraph contains exactly $2n$ edges, so in this case the geometric system
parameterizes the general leaf of the quantum algebra (\ref{sl2q}), (\ref{algebraMi-Mj}).
\end{remark}

\subsection{$R$-matrix representation for the quantum monodromy $A_n$-algebra}\label{ss:R-matrix}

The relations (\ref{sl2q}), (\ref{algebraMi-Mj}) can be presented in the $R$-matrix form. Assuming the standard
notation of the tensor product of matrix spaces 1 and 2, we have
\bea
&{}& R_{12}[q]\sheet{1}{M}_i R_{12}[q^{-1}]\sheet{2}{M}_j=\sheet{2}{M}_j R_{12}[q] \sheet{1}{M}_i R_{12}[q^{-1}],
\ i<j \label{R-ij}\\
&{}&
R^T_{12}[q^{-2}] \sheet{2}{M}_i \sheet{1}{M}_i= \sheet{1}{M}_i \sheet{2}{M}_i R_{12}[q^{-2}], \label{R-ii}
\eea
where
\be
R_{12}[q]:=q(E_{11}\otimes E_{11}+E_{22}\otimes E_{22})+E_{11}\otimes E_{22}+E_{22}\otimes E_{11}+
(q-q^{-1})E_{12}\otimes E_{21}
\label{R12}
\ee
is the standard $4\times 4$ quantum $R$-matrix by Kulish and Sklyanin (see \cite{KulSk}),
$$
R_{12}[q]=\left[
            \begin{array}{cccc}
              q &  &  &  \\
               & 1 & q{-}q^{-1} &  \\
               &  & 1 &  \\
               &  &  & q \\
            \end{array}
          \right].
$$
This matrix possesses the properties
\be
R_{12}[q^{-1}]=R^{-1}_{12}[q],\qquad R^T_{12}[q]=R_{21}[q],
\ee
where we let the superscript $T$ denote the total transposition.

Note that the reflection equation (\ref{R-ij}) results in the Jacobi identities for all distinct indices $i,j,k$
provided $R_{12}[q]$ satisfies the quantum Yang--Baxter equation
\be
R_{12}[q]R_{13}[q]R_{23}[q]=R_{23}[q]R_{13}[q]R_{12}[q].\label{YBE}
\ee
This follows from the chain of equalities for $i<j<k$ (for brevity, we omit the argument $q$ assuming $R_{12}:=R_{12}[q]$):
\bea
&{}&
[R_{12}\sheet1{M}_iR^{-1}_{12}\sheet2{M}_j]R^{-1}_{13}R^{-1}_{23}\sheet3{M}_k
=\sheet2{M}_jR_{12}\sheet1{M}_i[R^{-1}_{12}R^{-1}_{13}R^{-1}_{23}]\sheet3{M}_k
\nonumber\\
&{}&
=\sheet2{M}_jR_{12}R^{-1}_{23}[\sheet1{M}_iR^{-1}_{13}\sheet3{M}_k]R^{-1}_{12}
=\sheet2{M}_j[R_{12}R^{-1}_{23}R^{-1}_{13}]\sheet3{M}_k R_{13}\sheet1{M}_i R^{-1}_{13}R^{-1}_{12}
\nonumber\\
&{}&
=R^{-1}_{13}[\sheet2{M}_jR^{-1}_{23}\sheet3{M}_k]R_{12} R_{13}\sheet1{M}_i R^{-1}_{13}R^{-1}_{12}
=R^{-1}_{13}R^{-1}_{23}\sheet3{M}_k R_{23}\sheet2{M}_j[R^{-1}_{23}R_{12} R_{13}]\sheet1{M}_i R^{-1}_{13}R^{-1}_{12}
\nonumber\\
&{}&
=R^{-1}_{13}R^{-1}_{23}\sheet3{M}_k R_{23}R_{13}[\sheet2{M}_jR_{12}\sheet1{M}_i] [R^{-1}_{23}R^{-1}_{13}R^{-1}_{12}]
\nonumber\\
&{}&
=R^{-1}_{13}R^{-1}_{23}\sheet3{M}_k [R_{23}R_{13}R_{12}]\sheet1{M}_i R^{-1}_{12}
\sheet2{M}_jR_{12} R^{-1}_{12}R^{-1}_{13}R^{-1}_{23}
\nonumber\\
&{}&
=[R^{-1}_{13}R^{-1}_{23}R_{12}][\sheet3{M}_k R_{13}\sheet1{M}_i]R_{23} R^{-1}_{12}\sheet2{M}_j R^{-1}_{13}R^{-1}_{23}
\nonumber\\
&{}&
=R_{12}R^{-1}_{23}R^{-1}_{13}R_{13}\sheet1{M}_i R^{-1}_{13}[\sheet3{M}_k [R_{13}R_{23} R^{-1}_{12}]\sheet2{M}_j] R^{-1}_{13}R^{-1}_{23}
\nonumber\\
&{}&
=R_{12}R^{-1}_{23}\sheet1{M}_i [ R^{-1}_{13}R^{-1}_{12}R_{23}]\sheet2{M}_j R^{-1}_{23}\sheet3{M}_k
=R_{12}\sheet1{M}_iR^{-1}_{12}\sheet2{M}_j R^{-1}_{13}R^{-1}_{23}\sheet3{M}_k.
\nonumber
\eea
The proof of Jacobi identities for the triple product of entries of the same matrix $M_i$ follows from the
standard Yang--Baxter equation (\ref{YBE}) for the matrices $R_{12}[q^{-2}]$.

It seems however doubtful that we can write a similar $R$-matrix proof of Jacobi identities in the case
where we have two matrices $M_i$ and one $M_j$ with $i\ne j$. In fact if we did, then we would be able
to prove  the Jacobi identities for  arbitrary matrices $M_i$, which
is obviously impossible since no such closed algebra exists except a special case of the Painlev\'e VI algebra
considered in Sec.~\ref{ss:PVI-monodromy} below.

\subsection{Braid group action}\label{ss:braid}

As in \cite{ChM}, for $q=1$, we identify the matrices $M_i$ with the monodromy matrices of the Fuchsian
$2\times2$-matrix system:
\be
\partial_\lambda \Phi(\lambda)=\sum_{i=0}^{n-1}\frac{A_i}{\lambda-u_i} \Phi(\lambda);\quad \hbox{eigen\,}
A_i=\{1/4,-1/4\}.
\label{Dub-A}
\ee
Here $M_i$ are monodromy matrices of the 2-vector $\Phi(\lambda)$ when going around the $i$th singular point $u_i$.
Choosing $M_0=F_0$ corresponds to choosing a special basis in the space of solutions: in the
vicinity of every singular point $u_i$, we have two special solutions $\Phi_i^\pm(\lambda)=(\lambda-u_i)^{\pm1/4}f^\pm_{i}(\lambda)$, with $f^\pm_{i}(\lambda)$ being analytic in the
neighborhood of $u_i$. Then, obviously, on the classical level our choice of monodromy matrices corresponds
to choosing $\Phi^T(\lambda)=(\Phi_0^+(\lambda)+\Phi_0^-(\lambda),\Phi_0^+(\lambda)-\Phi_0^-(\lambda))$. The quantum
algebra (\ref{R-ij}), (\ref{R-ii}) is then the (consistent) quantum analogue of the Korotkin--Samtleben
Poisson algebra \cite{KS}.

The action of the braid group on monodromy matrices was proposed by Dubrovin and M.M. \cite{DM}. We have the
following lemma in the quantum case.
\begin{lemma}\label{braid}
The transformations $\beta_{i,i+1}$, $i=1,\dots,n-2$, such that
\be
\beta_{i,i+1}[M_j]=\left\{\begin{array}{ll}
                            M_j & j\ne i,i+1 \\
                            M_i & j=i+1 \\
                            -M_i M_{i+1}M_i & j=i
                          \end{array}\right.
                          \label{braid-M}
\ee
is an automorphism of the quantum algebra (\ref{R-ij}), (\ref{R-ii}) and satisfies the quantum braid-group
identities $\beta_{i,i+1}\beta_{i+1,i+2}\beta_{i,i+1}=\beta_{i+1,i+2}\beta_{i,i+1}\beta_{i+1,i+2}$.
We again assume the natural ordering of quantum operators in (\ref{braid-M}).
\end{lemma}

{\em Proof.} The proof of that transformations (\ref{braid-M}) satisfy the braid group relation repeats the
proof of these relations in the classical case (see \cite{DM}) provided we have the natural ordering of quantum
operators. The proof of that the transformation (\ref{braid-M}) is an automorphism of the quantum algebra
of monodromy matrices is the direct calculation in course of which we have found another nice representation of
the braid group action in terms of the quantum geodesic function (\ref{Gij}).

\vskip 1mm
The following results can be proved by straightforward computations:

\begin{prop}\label{prop-det}
The quantum determinant relations $b_ic_i-q^2a_i^2=1$, $i=1,\dots,n-1$, are preserved by the quantum
braid-group action (\ref{braid-M}).
\end{prop}

\begin{lemma}\label{lm:braid-Gii+1}
We can present the braid-group action (\ref{braid-M}) in the form
\be
\beta_{i,i+1}[M_j]=\left\{\begin{array}{ll}
                            M_j & j\ne i,i+1 \\
                            M_i & j=i+1 \\
                            q M_i G_{i,i+1}-q^2 M_{i+1}=q^{-1}G_{i,i+1}M_{i}-q^{-2}M_{i+1}& j=i,
                          \end{array}\right.
                          \label{braid-M-Gii+1}
\ee
and the quantum geodesic operator $G_{i,i+1}$ has the following commutation relations with
the monodromy matrices $M_j$:
\bea
&{}& q^{-1}G_{i,i+1}M_i-q M_iG_{i,i+1}=(q^{-2}-q^2)M_{i+1},\nonumber\\
&{}& q G_{i,i+1}M_{i+1}-q^{-1} M_{i+1}G_{i,i+1}=(q^{2}-q^{-2})M_{i},\label{GM}\\
&{}& G_{i,i+1}M_{j}- M_{j}G_{i,i+1}=0,\quad j\ne i,i+1.\nonumber
\eea
In general, the quantum geodesic functions $G_{i,j}$, $1\le i<j\le n-1$, have the following commutation
relations with the quantum monodromy matrices $M_k$, $k=1,\dots,n-1$:
\bea
&{}& q^{-1}G_{i,j}M_i-q M_iG_{i,j}=(q^{-2}-q^2)M_{j},\nonumber\\
&{}& q G_{i,j}M_{j}-q^{-1} M_{j}G_{i,j}=(q^{2}-q^{-2})M_{i},\label{GMij}\\
&{}& G_{i,j}M_{k}- M_{k}G_{i,j}=(q^{2}-q^{-2})[M_iG_{k,j}-G_{i,k}M_j],\quad i<k<j,\nonumber\\
&{}& G_{i,j}M_{k}- M_{k}G_{i,j}=0,\quad i<j<k\ \hbox{or}\ k<i<j. \nonumber
\eea
\end{lemma}

\begin{corollary}\label{cor:central}
The relations (\ref{GM}) together with the invertibility condition (\ref{M-square}) imply that
\be
G_{i,i+1}M_iM_{i+1}=M_iM_{i+1} G_{i,i+1}\quad\hbox{and}\quad G_{i,i+1}M_{i+1}M_i=M_{i+1} M_i G_{i,i+1}.
\label{GMM}
\ee
Therefore the product $M_1M_2\cdots M_{n-2}M_{n-1}$ and its inverse $M_{n-1}M_{n-2}\cdots M_2M_1$ are
the braid-group invariants.
\end{corollary}

\begin{remark}\label{rm:central}
The fact that we have the braid-group invariants (\ref{GMM}) besides the ``genuine'' Casimirs given by
Proposition~\ref{prop-det} is due to that we consider here only the subgroup of the total braid group that
leaves invariant the rooted edge.
\end{remark}

\subsection{Riemann sphere with four holes/orbifold points}\label{ss:PVI-monodromy}

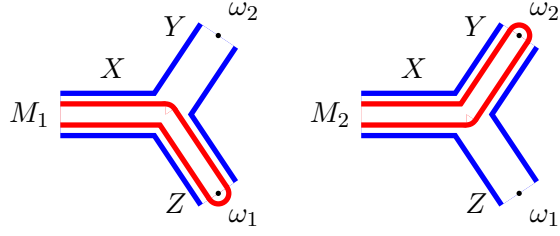
\begin{figure}[tb]
\begin{pspicture}(-4,-1.5)(4,1.5){
\newcommand{\MONE}{%
{\psset{unit=0.7}
\psline[linewidth=18pt,linecolor=blue](-2,0)(0,0)
\psline[linewidth=18pt,linecolor=blue](0,0)(1,1.5)
\psline[linewidth=18pt,linecolor=blue](0,0)(1,-1.5)
\psline[linewidth=14pt,linecolor=white](-2,0)(0,0)
\psline[linewidth=14pt,linecolor=white](0,0)(1,1.5)
\psline[linewidth=14pt,linecolor=white](0,0)(1,-1.5)
%
\rput(1,-1.5){\pscircle[linewidth=2pt,linecolor=red]{0.25}}
\rput(0,0){\pscircle[linewidth=2pt,linecolor=red]{0.25}}
\psline[linewidth=10pt,linecolor=red](-2,0)(0,0)
\psline[linewidth=10pt,linecolor=red](0,0)(1,-1.5)
\psline[linewidth=6pt,linecolor=white](-2,0)(0,0)
\psline[linewidth=6pt,linecolor=white](0,0)(1,-1.5)
\rput(1,1.5){\pscircle*{0.05}}
\rput(1,-1.5){\pscircle*{0.05}}
\rput(0.4,1.5){\makebox(0,0)[rb]{$Y$}}
\rput(-1,0.7){\makebox(0,0)[cb]{$X$}}
\rput(0.4,-1.5){\makebox(0,0)[rt]{$Z$}}
\rput(1.2,1.8){\makebox(0,0)[lb]{$\omega_2$}}
\rput(1.2,-1.8){\makebox(0,0)[lt]{$\omega_1$}}
\rput(-2.2,0){\makebox(0,0)[rc]{$M_1$}}
}
}
\newcommand{\MTWO}{%
{\psset{unit=0.7}
\psline[linewidth=18pt,linecolor=blue](-2,0)(0,0)
\psline[linewidth=18pt,linecolor=blue](0,0)(1,1.5)
\psline[linewidth=18pt,linecolor=blue](0,0)(1,-1.5)
\psline[linewidth=14pt,linecolor=white](-2,0)(0,0)
\psline[linewidth=14pt,linecolor=white](0,0)(1,1.5)
\psline[linewidth=14pt,linecolor=white](0,0)(1,-1.5)
%
\rput(1,1.5){\pscircle[linewidth=2pt,linecolor=red]{0.25}}
\rput(0,0){\pscircle[linewidth=2pt,linecolor=red]{0.25}}
\psline[linewidth=10pt,linecolor=red](-2,0)(0,0)
\psline[linewidth=10pt,linecolor=red](0,0)(1,1.5)
%
\psline[linewidth=6pt,linecolor=white](-2,0)(0,0)
\psline[linewidth=6pt,linecolor=white](0,0)(1,1.5)
%
\rput(1,1.5){\pscircle*{0.05}}
\rput(1,-1.5){\pscircle*{0.05}}
\rput(0.4,1.5){\makebox(0,0)[rb]{$Y$}}
\rput(-1,0.7){\makebox(0,0)[cb]{$X$}}
\rput(0.4,-1.5){\makebox(0,0)[rt]{$Z$}}
\rput(1.2,1.8){\makebox(0,0)[lb]{$\omega_2$}}
\rput(1.2,-1.8){\makebox(0,0)[lt]{$\omega_1$}}
\rput(-2.2,0){\makebox(0,0)[rc]{$M_2$}}
}
}
\rput(1.5,0){\MONE}
\rput(5.5,0){\MTWO}
}
\end{pspicture}
\caption{\small
The graphical representation for the monodromy matrices $M_1$ and $M_2$ in the
P$_{VI}$ case.}
\label{fi:M1-M2-PVI}
\end{figure}

In the case of a Riemann sphere with four holes/orbifold points, the fat graph has a  tree-like structure, where the corresponding orbifold points are now of
arbitrary orders; we therefore have arbitrary parameters $\omega_0$, $\omega_1$ and $\omega_2$. Fixing the zeroth hole/orbifold points as root, the arrangement is as in Fig.~\ref{fi:M1-M2-PVI} and we are left with three matrices: $F_{\omega_0}$ and another two, more complicated ones, which we call $M_1$ and $M_2$.
The quantum expressions for $M_1$ and $M_2$ read:
\be
M_1=X_XLX_ZF_{\omega_1}X_ZRX_X=\left[
                                   \begin{array}{cc}
                                     qa_1+\omega_1 & -b_1 \\
                                     c_1 & -q^{-1}a_1 \\
                                   \end{array}
                                 \right],
                                 \label{M1-PVI}
\ee
where $a_1=e^{-Z}$, $b_1=e^{X-Z}+e^{X+Z}+\omega_1 e^X$, and $c_1=e^{-X-Z}$, and
\be
M_2=X_XRX_YF_{\omega_2}X_YLX_X=\left[
                                   \begin{array}{cc}
                                     qa_2+\omega_2 & -b_2 \\
                                     c_2 & -q^{-1}a_2 \\
                                   \end{array}
                                 \right],
                                 \label{M2-PVI}
\ee
where $a_1=e^{Y}$, $b_1=e^{X+Y}$, and $c_1=e^{-X-Y}+e^{-X+Y}+\omega_2 e^{-X}$. The entries of both these
matrices satisfy the deformed $U_q(\mathfrak{sl}_2)$ relations:
\bea\label{defUq}
&{}&qa_ib_i=q^{-1}b_ia_i,\label{alg1}\\
&{}&q^{-1}a_ic_i=qc_ia_i,\qquad i=1,2 \label{sl2q-PVI}\\
&{}&b_ic_i=1+\omega_i q a_i+ q^2a_i^2,\quad c_ib_i=1+\omega_i q^{-1} a_i+ q^{-2}a_i^2.\label{alg3}
\eea

The algebraic relations between entries of $M_1$ and $M_2$ are
\bea
&{}& q^{-1}a_1a_2=qa_2a_1,\label{Mon-PVI1}\\
&{}& q^{-1}b_1b_2=qb_2b_1,\\
&{}& q^{-1}c_1c_2=qc_2c_1,\\
&{}& b_1a_2+q^{-2}a_1b_2+\omega_1 q^{-1}b_2 =a_2b_1 +q^2b_2a_1+\omega_1 q b_2,\\
&{}& a_1b_2=b_2a_1,\\
&{}& a_1c_2+q^{-2}c_1a_2+\omega_2 q^{-1}c_1 =c_2a_1+q^2a_2c_1+\omega_2 q c_1,\\
&{}& c_1a_2=a_2c_1,\\
&{}& qc_1b_2=q^{-1}b_2c_1,\\
&{}&qb_1c_2-q^{-1}c_2b_1=
[q^2-q^{-2}][q a_1a_2+q^{-1}a_2a_1+\omega_1a_2+\omega_2 a_1]+[q-q^{-1}]\omega_1\omega_2.\label{Mon-PVI9}
\eea
These relations in general are not consistent (otherwise we were be able to
formulate a monodromy algebra of $A_n$-type for arbitrary parameters $\omega_i$). They become
consistent if we supply them with an additional condition
\be
a_1a_2=q^2 c_1b_2.\label{addit-cond}
\ee
We call the algebraic relations (\ref{alg1})--(\ref{alg3}), (\ref{Mon-PVI1})--(\ref{Mon-PVI9}),
and (\ref{addit-cond}) the monodromy algebra for the P$_{VI}$ equation.

Note that we have three relations, (\ref{alg3}) for $i=1,2$ and (\ref{addit-cond}),
on six variable $a_i$, $b_i$, $c_i$,
which implies that our algebra must have an additional central element. It is not difficult to find that the elements
\bea
K_1&:=&a_1c_2-q^2c_1a_2-q\omega_2c_1\label{K1}\\
K_2&:=&a_2b_1-q^{-2}b_2a_1-q^{-1}\omega_1b_2\label{K2}
\eea
are central and they are related:
\be
K_1K_2=1.
\ee
We can thus identify the parameter $\omega_3$ corresponding to the fourth hole/orbifold point:
\be
\omega_3:=K_1+K_2.\label{omega3}
\ee

\begin{remark}\label{PVI-R}
It is interesting to mention that the monodromy matrices $M_1$ and $M_2$ given by the respective expressions
(\ref{M1-PVI}) and (\ref{M2-PVI}) satisfy the very same relation (\ref{R-ij}) with $i=1$ and $j=2$ and with
exactly the same $R$-matrix (\ref{R12}) provided the relations (\ref{Mon-PVI1})--(\ref{Mon-PVI9})
and (\ref{addit-cond}) be satisfied.
\end{remark}

We now define the three geodesic functions of the P$_{VI}$ system:
\bea
G_{XZ}&=&c_1+b_1+\omega_0 a_1,\label{GXZ}\\
G_{XY}&=&c_2+b_2+\omega_0 a_2,\label{GXY}\\
G_{YZ}&=&qb_1c_2-q^3a_1a_2-q^2(\omega_1a_2+\omega_2a_1)-q\omega_1\omega_2,\label{GYZ}
\eea
where we have introduced the last missing parameter $\omega_0$. For these quantities
we then obtain the system of algebraic relations of the P$_{VI}$ algebra of geodesic functions:
\bea\label{eq:AW(3)}
q G_{XY}G_{XZ}-q^{-1}G_{XZ}G_{XY}=[q^2-q^{-2}]G_{YZ}+[q-q^{-1}][\omega_1\omega_2+\omega_0\omega_3],\nonumber\\
q G_{XZ}G_{YZ}-q^{-1}G_{YZ}G_{XZ}=[q^2-q^{-2}]G_{XY}+[q-q^{-1}][\omega_2\omega_0+\omega_1\omega_3],\label{PVI}\\
q G_{YZ}G_{XY}-q^{-1}G_{XY}G_{YZ}=[q^2-q^{-2}]G_{XZ}+[q-q^{-1}][\omega_0\omega_1+\omega_2\omega_3].\nonumber
\eea

\begin{remark}
In \cite{MM}, it was proved that $F_{\omega_0}, M_1,M_2$ and a fourth matrix corresponding to the quantization of their product satisfy the Cherednik algebra $\mathcal H$ relations of type $\check{C_1}C_1$. The algebra (\ref{eq:AW(3)}) corresponds to the spherical subalgebra $e\mathcal He$  of
$\mathcal H$. In this paper we have shown that $M_1$ and $M_2$ satisfy the  deformed $U_q(\mathfrak{sl}_2)$ relations (\ref{defUq}), thus providing an interesting relation between the  Cherednik algebra of type $\check{C_1}C_1$ and the  deformed $U_q(\mathfrak{sl}_2)$. \end{remark}

\section{Conclusion}
In conclusion, we recall the main new results obtained in this paper.

We have quantized the MCG transformations (or flip morphisms) for flipping pending edges
of graphs corresponding to open Riemann surfaces with orbifold points.

Using the fact that elements of Fuchsian groups transform uniformly under the action of quantum MCG
transformations, we were able to provide
the explicit quantum ordering for a special class of quantum geodesic functions
corresponding to geodesics joining exactly two orbifold points or holes on a non-compact Riemann surface.

We have used the same quantum ordering for quantizing matrix elements of the Fuchsian group
associated to the Riemann surface. Interestingly, in all the cases where
it is possible to close the Poisson algebra of geodesic functions, such as
the  $A_n$-algebra related to
Schlesinger systems~\cite{Dub} and in the case of algebra related to the Painlev\'e VI equation~\cite{ChM-D4}, we have obtained a well defined quantum algebra for the elements in the corresponding Fuchsian group.

In the case of the  $A_n$-algebra, for each fixed geodesic, the matrix entries of the corresponding element in the Fuchsian group satisfy the quantum universal enveloping algebra $U_q(\mathfrak{sl}_2)$ relations, while in the case of algebra related to the Painlev\'e VI equation they satisfy a deformed version of $U_q(\mathfrak{sl}_2)$.

This result is quite interesting as it sheds light on the relation between the quantum universal enveloping algebra $U_q(\mathfrak{sl}_2)$ and the Zhedanov algebra $AW(3)$ already explored in \cite{WZ,Ter}. We have also obtained
the quantum generalisation of the Korotkin--Samtleben $r$-matrix Poisson algebra for the Schlesinger system
and constructed the action of the quantum braid group on the entries of quantum monodromy matrices.

In a separate publication, we shall consider the last remaining case of $D_n$-algebras of matrix elements
of the Fuchsian group for an annulus with $n$ ${\mathbb Z}_2$-orbifold points. The corresponding algebra of
quantum geodesic functions is related to the reflection equation with the spectral parameter.

\section*{Acknowledgments}
The authors are grateful to P. Wiegmann and A. Zabrodin for useful conversations. This work was supported by the Engineering and Physics Sciences Research Council EP/J007234/1. The work of L.Ch. was also supported in part by the Russian Foundation for Basic Research (Grant Nos. 11-01-00440-a and 13-01-12405-ofi-m), by the Grant of Supporting Leading Scientific Schools of the
Russian Federation NSh-4612.2012.1, and by the Program Mathematical Methods for Nonlinear
Dynamics.

\end{document}